\documentclass[reqno]{amsart}
\usepackage{bbm}
\usepackage{mathrsfs}
\usepackage{amsmath}
\allowdisplaybreaks[0]
\usepackage{paralist}
\usepackage{graphics}
\usepackage{epsfig}
\usepackage[pdftex]{hyperref}
\usepackage{amsfonts}
\usepackage{CJK}
\usepackage{fancyhdr}
\usepackage{graphicx}
\usepackage[dvipsnames,usenames]{color}
\usepackage{titletoc}
\usepackage{latexsym}
\usepackage{amssymb}
\usepackage{multicol}
\usepackage{graphics}
\usepackage{subfigure}
\usepackage{indentfirst}
\usepackage{cases}
\usepackage{curves}
\usepackage{cite}

\newtheorem{theorem}{Theorem}[section]

\newtheorem{lemma}{Lemma}[section]
\newtheorem{remark}{Remark}[section]
\newtheorem{example}[theorem]{Example}

\numberwithin{equation}{section}

\def\x#1{(\ref{#1})}
\def\R{{\Bbb R}}

\def\N{{\Bbb N}}

\def\E{{\mathbb{E}}}

\def\bc{\begin{center}}
\def\ec{\end{center}}
\def\ba{\begin{array}}
\def\ea{\end{array}}
\def\be{\begin{equation}}
\def\ee{\end{equation}}
\def\bea{\begin{eqnarray}}
\def\eea{\end{eqnarray}}
\def\beaa{\begin{eqnarray*}}
\def\eeaa{\end{eqnarray*}}

\def\ben{\begin{enumerate}}
\def\een{\end{enumerate}}
\def\hh{\!\!\!\!}
\def\EM{\hh &   &\hh}
\def\EQ{\hh & = & \hh}

\def\LE{\hh & \le & \hh}

\def\nn{\nonumber}
\def\oo{\infty}
\def\ifl{\iffalse}
\def\lb{\label}
\def\prf{\mbox{\bf Proof.~}}

\title[]{Robustness of nonuniform mean-square exponential dichotomies}

\author[Hailong Zhu]
{Hailong Zhu $^{1}$}
\address{$^1$ Anhui
University of Finance and Economics, Bengbu, 233030, China}

\email{hai-long-zhu@163.com (H. Zhu)}

\thanks{Hailong Zhu was supported by the National NSF of China (NO. 11671118), NSF of Anhui Province of China(NO. KJ2017A432, NO. KJ2018A0437).  }

\subjclass[2000]{60H10, 34D09} \keywords{Robustness; Nonuniform mean-square
exponential contraction; Nonuniform mean-square exponential dichotomy;
Stochastic differential equations.}

\begin{document}
\setlength{\parskip}{0.5\baselineskip}

\begin{abstract}
For linear stochastic differential equations (SDEs) with bounded coefficients, we
establish the robustness of  nonuniform mean-square exponential dichotomy (NMS-ED)
 on $[t_{0},+\oo)$, $(-\oo,t_{0}]$ and the  whole $\R$ separately, in the sense that
 such an NMS-ED persists  under a sufficiently small linear perturbation. The result for the
  nonuniform mean-square exponential contraction  (NMS-EC) is also discussed.
  Moreover,  in the process of   proving the existence of NMS-ED,
 we use the observation that the projections of the ``exponential growing solutions" and the ``exponential decaying solutions" on $[t_{0},+\oo)$, $(-\oo,t_{0}]$ and $\R$ are different but related. Thus, the relations of three types of projections on  $[t_{0},+\oo)$, $(-\oo,t_{0}]$ and $\R$ are discussed.

\end{abstract}

\maketitle

\section{\bf Introduction}
\setcounter{equation}{0} \noindent

The well-established notion of exponential dichotomy used in the analysis of nonautonomous systems is essentially   originated from the work of Perron \cite{per}. The theory of exponential dichotomy is   a powerful tool  to   describe  hyperbolicity   of   dynamical  systems  generated  by  differential   equations, especially for the stable and unstable invariant manifolds  of time-dependent systems. As mentioned in Coppel \cite{cop},``that dichotomies, rather than Lyapunov¡¯s characteristic exponents,
are the key to questions of asymptotic behaviour for nonautonomous differential equations".

Over   the   years, the classical exponential dichotomy and its properties have   been established  for evolution equations \cite{huy-06,  pa-88, pp-06, rr-95, lrs, ss}, functional differential   equations  \cite{cof-71,lin, per-71},
  skew-product   flows \cite{cl-94,cl,lmr,ss1} and random  systems  or   stochastic  equations \cite{drk,st, zlz, zc2, zcz}. We also refer to the books \cite{chi, cop, ms} for details and further references
related to exponential dichotomies.

However, dynamical  systems  exhibit  various  different   kinds  of   dichotomic    behavior and the classical notion of exponential dichotomy substantially restricts some dynamics.
In   order  to   investigate more general hyperbolicity, many attempts (see, e.g, \cite{np,np-97, pp-04})  have been made to extend the concept of classical  dichotomies.  Inspired by the work of  Barreira  and Pesin on the notion of nonuniformly  hyperbolic trajectory \cite{bp02,bp07}, Barreira and Valls  extended  the  concept  of   exponential   dichotomy to the nonuniform  ones and investigated some related problems, see for examples, the works \cite{bcv11, bcv13, bsv09,bv06,bv08} and the references therein.

On the other hand, from the point of view of  It\^{o} SDE, such properties of mean-square are  natural since the It\^{o} stochastic calculus is essentially deterministic in the mean-square setting, and there exist stationary coordinate changes under which flows of nonautonomous random differential equation can be viewed as those of SDE \cite{il-01}. Some related works on mean-square setting of random systems or stochastic equations can be found in \cite{fl,hig,hms,hmy,kl,ls,zc1}.
As our knowledge, mean-square exponential dichotomy
(MS-ED) was first introduced by Stanzhyts'kyi \cite{sk-01}, in which a sufficient condition has been proved to ensure that a linear SDE satisfies an MS-ED. Based on the definition of MS-ED, Stanzhyts'kyi and Krenevych \cite{sk-06}  proved the existence of a
quadratic form of linear SDE. In \cite{zc2} the robustness of MS-ED for a linear SDE was established. Stoica \cite{st} studied  stochastic
cocycles  in Hilbert spaces. Recently, Doan et al. \cite{drk} considered the MS-ED spectrum for random dynamical system.

Now we recall the definition of MS-ED. Consider the following linear $n$-dimensional It\^{o} stochastic system
\be\lb{a1}dx(t)=A(t)x(t)dt+G(t)x(t)d\omega(t), \quad t\in I,\ee
where $I$  is  either  the  half line $[t_{0},+\oo)$, $(-\oo,t_{0}]$ or the  whole $\R$, and $A(t)=(A_{ij}(t))_{n \times n}$, $G(t)=(G_{ij}(t))_{n \times n}$
are continuous functions with real entries. Eq. \x{a1} is said to possess an \emph{MS-ED} if there exists a linear projection $P(t): L^{2}(\Omega, \R^{n}) \rightarrow L^{2}(\Omega, \R^{n})$ such that
\be\lb{c1}\Phi(t)\Phi^{-1}(s)P(s)=P(t)\Phi(t)\Phi^{-1}(s),  \quad
\forall~ t,s\in I,\ee and positive constants $K, \alpha$  such that  \be\lb{a5}
\begin{split}
\E\|\Phi(t)\Phi^{-1}(s)P(s)\|^2 \leq K e^{-\alpha(t-s)},
\quad \forall~ (t,s)\in I^{2}_{\geq},\nn\\
\E\|\Phi(t)\Phi^{-1}(s)Q(s)\|^2 \leq K e^{-\alpha(s-t)}, \quad
\forall~ (t,s)\in I^{2}_{\leq},\nn
\end{split}
\ee
where  $\Phi(t)$ is a fundamental matrix solution of {\rm\x{a1}}, and $Q(t)={\rm Id}-P(t)$ is the
complementary projection of $P(t)$ for each $t\in I$.
$I^{2}_{\geq}:=\{(t,s)\in I^{2}: t\geq s\}$ and
$I^{2}_{\leq}:=\{(t,s)\in I^{2}: t\leq s\}$ denote the relations of $s$ and $t$ on $I$.

Inspired by the above,  this paper is to study the robustness of NMS-ED.  \x{a1} is said to possess an  \emph{NMS-ED} if  there exist a linear projection $P(t): L^{2}(\Omega, \R^{n}) \rightarrow L^{2}(\Omega, \R^{n})$  such that \x{c1} holds, and some  constants $M, \alpha>0$, $\varepsilon \ge 0$ such that
\begin{eqnarray}
&&\E\|\Phi(t)\Phi^{-1}(s)P(s)\|^2 \leq M e^{-\alpha(t-s)+\varepsilon |s|},
\quad \forall~ (t,s)\in I^{2}_{\geq},
\lb{c2}
\\
&&
\E\|\Phi(t)\Phi^{-1}(s)Q(s)\|^2 \leq M e^{-\alpha(s-t)+\varepsilon |s|}, \quad
\forall~ (t,s)\in I^{2}_{\leq},
\lb{x2}
\end{eqnarray}
where $\Phi(t)$ is a fundamental matrix solution of {\rm\x{a1}},  $Q(t)={\rm Id}-P(t)$ is the
complementary projection of $P(t)$ for each $t\in I$.
$I^{2}_{\geq}:=\{(t,s)\in I^{2}: t\geq s\}$ and
$I^{2}_{\leq}:=\{(t,s)\in I^{2}: t\leq s\}$ denote the relations of $s$ and $t$ on $I$.
For convenience, the  constants $\alpha$ and  $K$ in \x{c2}-\x{x2} are called
the   \emph{exponent} and   the   \emph{bound} of the NMS-ED respectively, as  in  the case of  deterministic systems \cite{hen-81}.
$\varepsilon$ is called the \emph{nonuniform degree} of the NMS-ED. In particular,
while $\varepsilon=0$, we obtain the notion of (uniform) MS-ED. We refer to \cite{sk-01,sk-06,st,zc1,zc2,zcz} for related results and techniques about this topic.

It is clear
that the notion of NMS-ED is a weaker requirement in comparison to the notion of MS-ED.
Actually, there exists a linear SDE which has an NMS-ED with nonuniform degree $\varepsilon$ cannot be removed.
For example, let $a>b>0$ be real parameters,
  \[
\left\{ \begin{array}{ll}
du & =(-a-bt\sin t)u(t)dt+\sqrt{2b\cos t}\exp(-at+bt\cos t)d\omega(t),\\
dv & =(a+bt\sin t)v(t)dt-\sqrt{2b\cos t}\exp(at-bt\cos t)d\omega(t)
\end{array} \right.
\]
admits an NMS-ED  which is not uniform. See Example \ref{exp51} in Section 6 for details.

Robustness (also known as roughness , see, e.g., \cite{cop}) here means that an NMS-ED persists  under a sufficiently small linear perturbation.
More precisely, for small  perturbations $B$, $H$, the following linear SDE
\be\lb{a3}dy(t)=(A(t)+B(t))y(t)dt+(G(t)+H(t))y(t)d\omega(t)\ee
also   admits  an NMS-ED.  As indicated by Coppel (\cite[p. 28]{cop}),  the   robustness  of
exponential  dichotomies was first   proved
by Massera and  Sch\"{a}ffer \cite{ms}, which states  that  all  ``neighboring"   linear systems  also   have   the same   dichotomy  with   a   similar   projection if the same happens for the original system. Robustness is one of the most basic  concepts  appearing in the theoretical studies of dynamical systems. This topic plays a key role in the stability theory for dynamical systems. For some early papers  about robustness (with the exception of  \cite{cop} and \cite{ms} mentioned above) are due to  Dalec'ki\u\i ~ and  Kre\u\i n \cite{dk-74}, and Palmer \cite{pa-84} for ordinary differential equations, Henry \cite{hen-81}, and Lin \cite{lin-94} for parabolic partial differential equations, Hale and Lin \cite{hl-86}, and Lizana \cite{liz-92} for functional differential equations, Pliss and Sell \cite{ps-99}, Chow and Leiva \cite {cl} for skew-product semiflow. For more recent works we mention in particular the papers \cite{bsv09, bv08, jw-01, pop-06, pop-09, zlz-17, zz-16}. It is worth mentioning  that on half line $\R^{+}$, $\R^{-}$ as well as the whole $\R$,
Ju and Wiggins \cite{jw-01}, and Popescu \cite{pop-06,pop-09} considered the case of roughness for exponential dichotomy  and analyze their dynamical behavior; Zhou, Lu, and Zhang \cite{zlz-17} discussed the relationship between nonuniform exponential
dichotomy and admissibility.

In this study, we extend the results and improve the method of \cite{zc2}. The main
differences of our results and those of \cite{zc2} are as follows:

\begin{itemize}
  \item In contrast to  \cite{zc2}, we extend the case of robustness of MS-ED to the general nonuniform setting. For this purpose, we need to  pass from small bounded perturbations of the coefficient matrix to exponentially decaying perturbations.
      \vspace{0.1cm}
  \item In \cite{zc2}, we only consider the case of robustness on the whole line $\R$. In the present paper, we prove the robustness of \x{a3} on half line $[t_{0},+\oo)$, $(-\oo,t_{0}]$ and the  whole $\R$. The proof is much more delicate than that of  MS-ED \cite{zc2}. This is because in different intervals, the different but related  explicit expressions  of the projections of the ``exponential growing solutions" and the ``exponential decaying solutions" for the perturbed equation \x{a3}  need first to be determined.
      \vspace{0.1cm}
  \item Furthermore, in contrast to paper \cite{zc2}, we analyze and compare the results obtained from operators that make up the projections of \x{a1} and \x{a3} on different intervals (see Theorem \ref{main32} and Remark \ref{rem51}), and estimate the distance between the solution of \x{a1} and the perturbed solution of \x{a3} (see Theorem \ref{main33} and Remark \ref{rem32}).
\end{itemize}

The rest part of this paper is organized as follows. The robustness of NMS-EC is established in Section 2.  Section 3  proves the robustness of
 NMS-ED on half line  $[t_{0},+\oo)$ and analyze that the solution of \x{a1} and the perturbed solution of \x{a3} are forward  asymptotic in the mean-square sense.
The robustness under the nonuniform setting on half line $(-\oo,t_{0}]$ is  presented in Section 4.   Section 5 combines the advantages of the projections on half line  $[t_{0},+\oo)$ and  $(-\oo,t_{0}]$, and proves the robustness of NMS-ED on the whole  $\R$. In addition, the relationship   of the projections on $[t_{0},+\oo)$, $(-\oo,t_{0}]$ and $\R$ is also discussed in Section 5.
Finally, an example  is given in Section 6,
which indicates that there exists a linear SDE which admits an NMS-ED but not uniform.

\section{Robustness of NMS-EC}
\setcounter{equation}{0} \noindent

In this section we will answer the following question: Does \x{a3} admit an NMS-EC if \x{a1} admits an NMS-EC while $B,~H$ is small?
That is to say, we consider the  robustness of NMS-EC. The following statement is a particular case of  NMS-ED with projection $P(t)=Id$ for every $t\in I$.
{\rm \x{a1}} is said to admit an \emph{NMS-EC} if for some  constants $M, \alpha>0$ and $\varepsilon \ge 0$ such that  \be\lb{b2}
\E\|\Phi(t)\Phi^{-1}(s)\|^2 \leq M e^{-\alpha(t-s)+\varepsilon |s|},
\quad \forall~ (t,s)\in I^{2}_{\geq}.
\ee
 In particular, when $\varepsilon=0$ in  \x{b2}, we obtain the notion of uniform mean-square exponential contraction.

Throughout this paper, we assume that $(\Omega, \mathscr{F},
\mathbb{P})$ is a  probability space,
$\omega(t)=(\omega_{1}(t),\ldots \omega_{n}(t))^{T}$ is an
$n$-dimensional Brownian motion defined on the space $(\Omega,
\mathscr{F}, \mathbb{P})$. $\|\cdot\|$ is used to denote both  the Euclidean vector norm or the matrix norm as appropriate, and $L^{2}(\Omega, \R^{n})$
stands for the space of all $\R^{n}$-valued random variables
$x: \Omega \rightarrow \R^{n}$ such that
$$\E \|x\|^2=\int_{\Omega}\|x\|^2 d \mathbb{P}<\oo.$$

In order to describe the robustness in an explicit form, we present the following theorem, which shows that the NMS-EC is robust under sufficiently
small linear perturbations. Here we mention that the NMS-EC considered in this section is in an arbitrary interval $I\subset \R$.

\begin{theorem}\label{main21} Let $A(\cdot), B(\cdot), G(\cdot), H(\cdot)$ be $n \times n$-matrix continuous functions with real entries such that
{\rm\x{a1}} admits an NMS-EC {\rm\x{b2}} with coefficient matrix bounded and perturbation exponential decaying in $I$,  i.e., there exist
constants $a, b, g, h>0$ such that
\beaa\lb{b4}\|A(t)\|\le a,\quad \|G(t)\|\le g,\quad \|B(t)\|\le b e^{-\frac{\varepsilon |t|}{2}}, \quad \|H(t)\|\le
h e^{-\frac{\varepsilon |t|}{2}},\quad t\in I.\eeaa
Let $b, h$ small enough such that
\bea\lb{b6}\tilde{M}:=8b^{2}+8g^{2}h^{2}+\alpha h^2 <\frac{\alpha^{2}}{6M}.\eea
 Then
{\rm\x{a3}} also  admits an NMS-EC in $I$ with the
bound $M$ replaced by $3M$, and exponent $\alpha$ replaced by $-\frac{\alpha}{2}+\frac{3M\tilde{M}}{\alpha}$, i.e.,
 \be\lb{b5}\E\|\hat{\Phi}(t)\hat{\Phi}^{-1}(s)\|^2\leq 3M
e^{(-\frac{\alpha}{2}+\frac{3M\tilde{M}}{\alpha})(t-s)+\varepsilon |s|}, \quad \forall~ (t,s)\in
I^{2}_{\geq},\ee where $\hat{\Phi}(t)$ is a fundamental matrix solution of {\rm\x{a3}}.
\end{theorem}

\noindent\prf{Write \[\hat{\Phi}(t,s)=\hat{\Phi}(t)\hat{\Phi}^{-1}(s).\] One can easily verify that $\hat{\Phi}(t,s)$  is a fundamental matrix solution of {\rm\x{a3}} with $\hat{\Phi}(s,s)=Id$.  $L^{2}(\Omega, \R^{n})$ is a Banach space with the norm $(\E \|x\|^2)^{\frac{1}{2}}$. The
Banach algebra of bounded linear operators on $L^{2}(\Omega, \R^{n})$ is denoted by $\mathfrak{B}(L^{2}(\Omega, \R^{n}))$.
Now we introduce the space
\be\lb{b7}\mathscr{L}_{c}:=\{\hat{\Phi}: I^{2}_{\geq} \rightarrow
\mathfrak{B}(L^{2}(\Omega, \R^{n})): ~\hat{\Phi} {\rm~is ~continuous ~and~}
\|\hat{\Phi}\|_{c}<\oo\}\ee with the norm \be\lb{b8}
\|\hat{\Phi}\|_{c}=\sup\left\{(\E\|\hat{\Phi}(t, s)\|^2)^{\frac{1}{2}}
e^{-\frac{\varepsilon}{2} |s|}:
(t,s)\in I^{2}_{\geq}\right\}.\ee
 Clearly, $(\mathscr{L}_{c},
\|\cdot\|_{c})$ is a Banach spaces. In order to state our result, we need the following existence and uniqueness lemma.
\begin{lemma}\lb{lem23}
  For  any given initial value $\xi_{0} \in \R^{n}$, {\rm \x{a3}} has a unique solution $\hat{\Phi}(t,s) \xi_{0}$ with $\hat{\Phi} \in (\mathscr{L}_{c},
\|\cdot\|_{c})$ such that
   \bea\lb{b9} \hat{\Phi}(t,s)\EQ
\Phi(t)\Phi^{-1}(s)+ \int^{t}_{s}\Phi(t)\Phi^{-1}(\tau)H(\tau)\hat{\Phi}(\tau, s)d\omega(\tau)
\nn\\
\EM + \int^{t}_{s}\Phi(t)\Phi^{-1}(\tau)\big{(}B(\tau)-
G(\tau)H(\tau)\big{)}\hat{\Phi}(\tau, s)d\tau\eea
   with $ \hat{\Phi}(s,s) \xi_{0}=\Phi(s)\Phi^{-1}(s) \xi_{0}= \xi_{0}$.
\end{lemma}

\prf{ In what follows (in order to simplify the presentation), write $\tilde{B}(t)=B(t)-G(t)H(t)$.
We first prove that the function $\hat{\Phi}(t,s) \xi_{0}$  is a solution of \x{a3}.
Set
\beaa \xi(t)\EQ
\Phi^{-1}(s)\xi_{0}+ \int^{t}_{s}\Phi^{-1}(\tau)H(\tau)\hat{\Phi}(\tau, s)\xi_{0}d\omega(\tau)
\\
\EM + \int^{t}_{s}\Phi^{-1}(\tau)\tilde{B}(\tau)\hat{\Phi}(\tau, s)\xi_{0}d\tau.\eeaa
Let $y(t)=\Phi(t)\xi(t)$. Clearly,
\[ \hat{\Phi}(t,s) \xi_{0}= \Phi(t)\xi(t) =y(t).
\]
One can easily verify that $\xi(t)$ satisfies the differential
\beaa d \xi(t)\EQ \Phi^{-1}(t)\big{(}B(t)-
G(t)H(t)\big{)}y(t)d t +\Phi^{-1}(t)H(t)y(t)d\omega(t).
\eeaa
Since $\Phi(t)$ is a fundamental matrix solution of \x{a1}, it follows from It\^{o} product rule that
\beaa dy(t)\EQ
d\Phi(t)\xi(t)+\Phi(t)d\xi(t)+G(t)\Phi(t)\Phi^{-1}(t)H(t)y(t)dt\\
\EQ A(t)y(t)dt+G(t)y(t)d\omega(t)+\big{(}B(t)-
G(t)H(t)\big{)}y(t)dt\\
\EM+H(t)y(t)d\omega(t)+G(t)H(t)y(t)dt\\\EQ
(A(t)+B(t))y(t)dt+(G(t)+H(t))y(t)d\omega(t),
\eeaa
which means that $y(t)=\hat{\Phi}(t,s)  \xi_{0}$ is a solution of \x{a3}.  This conclusion is consistent with that in \cite[Theorem 3.3.1]{mao} (see also \cite[Section 2.4.2]{lad}).

Now we prove that $\hat{\Phi}$ is unique in $(\mathscr{L}_{c},
\|\cdot\|_{c})$.  Let
\beaa(\Gamma \hat{\Phi})(t,s)\EQ
\Phi(t)\Phi^{-1}(s)+ \int^{t}_{s}\Phi(t)\Phi^{-1}(\tau)H(\tau)\hat{\Phi}(\tau, s)d\omega(\tau)
\nn\\
\EM + \int^{t}_{s}\Phi(t)\Phi^{-1}(\tau)\tilde{B}(\tau)\hat{\Phi}(\tau, s)d\tau.\eeaa
It follows from \x{b2}, $\E\|x\|\le \sqrt{\E\|x\|^{2}}$, Cauchy-Schwarz inequality, It\^{o} isometry property of stochastic integrals, and the
elementary inequality \be\lb{b10}\left\|\sum_{k=1}^{m}a_{k}\right\|^{2}\le
m\sum_{k=1}^{m}\|a_{k}\|^2\ee that
\bea \EM \E\|(\Gamma \hat{\Phi})(t,s)\|^{2} \le 3\E \|\Phi(t)\Phi^{-1}(s)\|^{2} +
3\E \left\|\int^{t}_{s}\Phi(t)\Phi^{-1}(\tau)H(\tau)\hat{\Phi}(\tau, s)d\omega(\tau)\right\|^{2}\nn\\
\EM + 3\E \left\| \int^{t}_{s}\Phi(t)\Phi^{-1}(\tau)\tilde{B}(\tau)\hat{\Phi}(\tau, s)d\tau \right\|^{2}\nn\\
\LE 3M e^{-a(t-s)+\varepsilon |s|}+3\int^{t}_{s} \E\|\Phi(t)\Phi^{-1}(\tau)\|^{2}
\E\|H(\tau)\|^{2} \E\|\hat{\Phi}(\tau, s)\|^{2} d\tau \nn\\
\EM + 3 \left(\int^{t}_{s} \E \|\Phi(t)\Phi^{-1}(\tau)\|\E\|\tilde{B}(\tau)\|d \tau\right) \times \left(\int^{t}_{s} \E \|\Phi(t)\Phi^{-1}(\tau)\| \E\|\tilde{B}(\tau)\| \E\|\hat{\Phi}(\tau, s)\|^{2}d \tau\right)\nn\\
\LE 3M e^{-\alpha(t-s)+\varepsilon |s|}+ 3Me^{\varepsilon |s|} \sup_{(\tau,s)\in I^{2}_{\ge}}\left(\E\|\hat{\Phi}(\tau, s)\|^{2}e^{-\varepsilon |s|}\right)\nn\\
 \EM
\times  \left(h^{2}\int^{t}_{s} e^{-\alpha(t-\tau)}d \tau + 2(b^{2}+g^{2}h^{2}) \left(\int^{t}_{s} e^{-\frac{\alpha}{2}(t-\tau)}d \tau\right)^{2}\right)\nn\\
\LE 3M e^{\varepsilon |s|} +3M e^{\varepsilon |s|} (\frac{\alpha h^{2}+8b^{2}+8g^{2}h^{2}}{\alpha^{2}})\sup_{(\tau,s)\in I^{2}_{\ge}}\left(\E\|\hat{\Phi}(\tau, s)\|^{2}e^{-\varepsilon |s|}\right),\nn
\eea
and this implies that
\[\E\|(\Gamma \hat{\Phi})(t,s)\|^{2}e^{-\varepsilon s}\le 3M +\frac{3\tilde{M}M}{\alpha^{2}}\sup_{(\tau,s)\in I^{2}_{\ge}}\left(\E\|\hat{\Phi}(\tau, s)\|^{2}e^{-\varepsilon |s|}\right)< \oo\]
with $\tilde{M}=8b^{2}+8g^{2}h^{2}+\alpha h^2$.
Proceeding in the
same procedure above, for any $\hat{\Phi}_{1}, \hat{\Phi}_{2}\in
\mathscr{L}_{c}$, we have
\be \lb{b12}\|\Gamma\hat{\Phi}_{1}-\Gamma\hat{\Phi}_{2}\|^{2}_{c}\le \frac{3\tilde{M}M}{\alpha^{2}}\sup_{(\tau,s)\in I^{2}_{\ge}}\left(\E\|\hat{\Phi}_{1}(\tau, s)-\hat{\Phi}_{1}(\tau, s)\|^{2}e^{-\varepsilon |s|}\right).\ee
Note that
\beaa\sup_{(t,s)\in I^{2}_{\ge}}\left(\E\|\hat{\Phi}_{1}(t, s)-\hat{\Phi}_{1}(t, s)\|^{2}e^{-\varepsilon |s|}\right)\EQ
\sup_{(t,s)\in I^{2}_{\ge}}\left((\E\|\hat{\Phi}_{1}(t, s)-\hat{\Phi}_{1}(t, s)\|^{2})^{\frac{1}{2}}e^{-\frac{\varepsilon |s|}{2} }\right)^{2}\\
\LE \left(\sup_{(t,s)\in I^{2}_{\ge}}(\E\|\hat{\Phi}_{1}(t, s)-\hat{\Phi}_{1}(t, s)\|^{2})^{\frac{1}{2}}e^{-\frac{\varepsilon |s|}{2} }\right)^{2}\\
\EQ \|\hat{\Phi}_{1}-\hat{\Phi}_{2}\|^{2}_{c},
\eeaa
which together with \x{b12} implies
\[\|\Gamma\hat{\Phi}_{1}-\Gamma\hat{\Phi}_{2}\|_{c}\le \sqrt{\frac{3\tilde{M}M}{\alpha^{2}} } \|\hat{\Phi}_{1}-\hat{\Phi}_{2}\|_{c}.\]
Since $\tilde{M}<\frac{\alpha^{2}}{3M}$, $\Gamma$ is a contraction operator. Hence, there exist a
unique $\hat{\Phi} \in \mathscr{L}_{c}$ such that $\Gamma\hat{\Phi}= \hat{\Phi}$, which satisfies the identity \x{b9}. This completes the proof of the lemma. \hspace{\stretch{1}}$\Box$}

We proceed with the proof of the theorem.  Squaring both sides of \x{b9}, and taking expectations, it follows from \x{b10} that
\bea\lb{b13} \E\|\hat{\Phi}(t,s)\|^{2} \LE
3\E\|\Phi(t)\Phi^{-1}(s)\|^{2}+ 3\E\left\|\int^{t}_{s}\Phi(t)\Phi^{-1}(\tau)H(\tau)\hat{\Phi}(\tau,s)
d\omega(\tau)\right\|^{2}\nn\\
\EM + 3\E\left\|\int^{t}_{s}\Phi(t)\Phi^{-1}(\tau)\tilde{B}(\tau)\hat{\Phi}(\tau,s)d\tau\right\|^{2}.
\eea
By using It\^{o} isometry
property and inequalities \x{b2}, the second term of right-hand side in  \x{b13} can be deduced as follows:
\beaa \EM \E\left\|\int^{t}_{s}\Phi(t)\Phi^{-1}(\tau)H(\tau)\hat{\Phi}(\tau,s)
d\omega(\tau)\right\|^{2} \\
\EQ \int^{t}_{s}\E\|\Phi(t)\Phi^{-1}(\tau)\|^{2} \E\|H(\tau)\|^{2} \E\|\hat{\Phi}(\tau,s)\|^{2} d \tau\\
\LE Mh^{2}\int^{t}_{s}e^{-\alpha(t-\tau)}\E\|\hat{\Phi}(\tau,s)\|^{2} d \tau.
\eeaa
As to the third term in \x{b13}, it follows from $\E\|x\|\le \sqrt{\E\|x\|^{2}}$, Cauchy-Schwarz inequality, and the inequalities \x{b2} that
\beaa \EM \E\left\|\int^{t}_{s}\Phi(t)\Phi^{-1}(\tau)\tilde{B}(\tau)\hat{\Phi}(\tau,s)d\tau\right\|^{2}\\
\EQ \E\left\|\int^{t}_{s}\left(\Phi(t)\Phi^{-1}(\tau)\tilde{B}(\tau)\right)^{\frac{1}{2}}
\left(\left(\Phi(t)\Phi^{-1}(\tau)\tilde{B}(\tau)\right)^{\frac{1}{2}}\hat{\Phi}(\tau,s)\right)
d\tau\right\|^{2}\\
\LE \left(\int^{t}_{s}\E\left\|\Phi(t)\Phi^{-1}(\tau)\tilde{B}(\tau)\right\| d\tau \right)\\
\EM \times \left(\int^{t}_{s}\E\left\|\Phi(t)\Phi^{-1}(\tau)\tilde{B}(\tau)\right\| \E\left\|\hat{\Phi}(\tau,s)\right\|^{2}d\tau \right)\\
\LE 2M(b^{2}+g^{2}h^{2})\left(\int^{t}_{s}e^{-\frac{\alpha}{2}(t-\tau)} d\tau \right)
\left(\int^{t}_{s}e^{-\frac{\alpha}{2}(t-\tau)} \E\left\|\hat{\Phi}(\tau,s)\right\|^{2}d\tau \right)\\
\LE  \frac{4M(b^{2}+g^{2}h^{2})}{\alpha}\int^{t}_{s}e^{-\frac{\alpha}{2}(t-\tau)} \E\left\|\hat{\Phi}(\tau,s)\right\|^{2}d\tau.
\eeaa
Since $\alpha>0$, we can rewrite the inequality \x{b13} as
\bea\lb{b14} \EM\E\|\hat{\Phi}(t,s)\|^{2} \le 3M e^{-\alpha(t-s)+ \varepsilon |s|} + 3Mh^{2}\int^{t}_{s}e^{-\alpha(t-\tau)}\E\left\|\hat{\Phi}(\tau,s)\right\|^{2} d \tau\nn\\
\EM +
\frac{12M(b^{2}+g^{2}h^{2})}{\alpha}\int^{t}_{s}e^{-\frac{\alpha}{2}(t-\tau)} \E\left\|\hat{\Phi}(\tau,s)\right\|^{2}d\tau\nn\\
\LE 3M e^{-\frac{\alpha}{2}(t-s)+ \varepsilon |s|} +3M (\frac{\alpha h^{2}+8b^{2}+8g^{2}h^{2}}{\alpha})
\int^{t}_{s}e^{-\frac{\alpha}{2}(t-\tau)} \E\left\|\hat{\Phi}(\tau,s)\right\|^{2}d\tau.\nn\\
\eea
%
Let
\be\lb{b15} x(t)=\E\|\hat{\Phi}(t,s)\|^{2},\quad X(t)=3M e^{-\frac{\alpha}{2}(t-s)+ \varepsilon |s|} + \frac{3M \tilde{M}}{\alpha} \int^{t}_{s} e^{-\frac{\alpha}{2}(t-\tau)} x(\tau) d\tau
\ee
for any fixed $s\in I$ with $\tilde{M}=\alpha h^{2}+8b^{2}+8g^{2}h^{2}$. Clearly, inequality \x{b14} can be rewritten as
\[x(t)\le X(t),\quad \mbox{for~all}~(t,s)\in
I^{2}_{\geq}.\]
On the other hand,
\[\frac{d}{dt}X(t)=-\frac{\alpha}{2}X(t)+\frac{3M \tilde{M}}{\alpha}x(t),\] and therefore, \[\frac{d}{dt}X(t)\leq \left(\frac{3M \tilde{M}}{\alpha}-\frac{\alpha}{2}\right) X(t).\]
Integrating the above inequality from $s$ to $t$ and note that
$X(s)=3Me^{\varepsilon |s|}$, we obtain \be\lb{b18}x(t)\leq X(t)\leq 3Me^{\varepsilon |s|}
e^{(-\frac{\alpha}{2}+\frac{3M\tilde{M}}{\alpha})(t-s)},\quad \mbox{for ~all}~ (t,s)\in
I^{2}_{\geq}.\ee By \x{b18}, using \x{b15}, we obtain the desired inequality \x{b5}, and
this completes the proof of the theorem. \hspace{\stretch{1}}$\Box$}

\begin{remark}
  Since the nonuniform degree $\varepsilon>0$ exists for $(t,s)\in
I^{2}_{\geq}$, the perturbations $B$ and $H$ should be chosen with exponential decaying to eliminate the effect caused by the nonuniform degree. For the uniform case, it suffices to consider the bounded condition instead of exponential decaying. See {\rm\cite{zc2}} for details about the case of $\varepsilon=0$, which generalizes (and imitates) the notion of robustness of exponential dichotomy for ODE (see e.g., {\rm\cite{cop, ms}}).
\end{remark}

As a special case of \x{a3}, if we
consider the system \be\lb{b16}dy(t)=(A(t)+B(t))y(t)dt+G(t)y(t)d\omega(t),\ee in
which the linear perturbed term only appears in the ``drift". Of course,
Theorem \ref{main21} can also be applied to \x{b16}  but merely with the development of slightly better estimation (with the bound and the exponent replaced by smaller constants) than the one in Theorem \ref{main21}, since there is no perturbation in the ``volatility". Actually, for  any given initial value $\xi_{0} \in \R^{n}$, {\rm \x{b16}} has a unique solution $\hat{\Phi}(t,s) \xi_{0}$ with $\hat{\Phi} \in (\mathscr{L}_{c},
\|\cdot\|_{c})$ such that
   \beaa\lb{b17} \hat{\Phi}(t,s)\EQ
\Phi(t)\Phi^{-1}(s)+  \int^{t}_{s}\Phi(t)\Phi^{-1}(\tau)B(\tau)\hat{\Phi}(\tau, s)d\tau\eeaa
 instead of \x{b9}, which is  more similar to solutions of the classical ordinary differential equations. See e.g., \cite{hale}.

\begin{theorem}\label{main22} Let $A(\cdot), B(\cdot), G(\cdot)$ be $n \times n$-matrix continuous functions with real entries such that
{\rm\x{a1}} admits an NMS-EC {\rm\x{b2}} with coefficient matrix bounded and perturbation exponential decaying in $I$,  i.e., there exist
constants $a, b, g>0$ such that
\[\|A(t)\|\le a,\quad \|G(t)\|\le g,\quad \|B(t)\|\le b e^{-\frac{\varepsilon |t|}{2}}, \quad t\in I.\]
If $b<\alpha/(2\sqrt{2M})$, then {\rm\x{b16}} also admits an NMS-EC in $I$ with the
bound $M$ replaced by $2M$, and exponent $\alpha$ replaced by $-\frac{\alpha}{2}+\frac{4Mb^{2}}{\alpha}$, i.e.,
 \[\E\|\hat{\Phi}(t)\hat{\Phi}^{-1}(s)\|^2\leq 2M
e^{(-\frac{\alpha}{2}+\frac{4Mb^{2}}{\alpha})(t-s)+\varepsilon |s|}, \quad \forall~ (t,s)\in
I^{2}_{\geq}.\]
\end{theorem}

\section{Robustness of NMS-ED on the half line $[t_{0},+\oo)$}
\setcounter{equation}{0} \noindent

In   this   section   we  state and prove our  main  result   on  the robustness of NMS-ED on $I=[t_{0},+\oo)$. The case of the interval $I=(-\oo,t_{0}]$ and the whole $\R$ will be discussed in Section 4 and Section 5 respectively.

The following theorem is on the robustness of NMS-ED of \x{a1} on $[t_{0},+\oo)$,
 and its proof is more general and complicated than that of Theorem \ref{main21}, because we need to find out the explicit expressions of the ``exponential growing solutions" and the ``exponential decaying solutions" for the perturbed equation \x{a3}
 along the stable and unstable directions respectively.  To do this, we  rewrite the unique solution of \x{a3} along the stable direction under a natural condition: boundedness.  It is also worth mentioning that the following theorem is also valid for NMS-EC. Indeed a contraction is a dichotomy with $P(t)=Id$ for every $t\in I$.

\begin{theorem}\label{main23} Let $A(\cdot), B(\cdot), G(\cdot), H(\cdot)$ be $n \times n$-matrix continuous functions with real entries such that
{\rm\x{a1}} admits an NMS-ED {\rm\x{c2}-\x{x2}} with $\varepsilon<\alpha$, and assume that coefficient matrices of {\rm\x{a3}} satisfy
\be\lb{x1}\|A(t)\|\le a,\quad \|G(t)\|\le g,\quad \|B(t)\|\le b e^{-\varepsilon |t|}, \quad \|H(t)\|\le
h e^{-\varepsilon |t|},\quad t\in I\ee
with constants $a, b, g, h>0$.
Let $b, h$ small enough such that
\beaa\lb{b6}\tilde{M}:=8b^{2}+8g^{2}h^{2}+\alpha h^2 <\frac{\alpha^{2}}{20M}.\eeaa
Then
{\rm\x{a3}} admits an NMS-ED in $I$ with
linear projections $\hat{P}(t): L^{2}(\Omega, \R^{n}) \rightarrow L^{2}(\Omega, \R^{n})$ such that
\be\lb{x4}\hat{\Phi}(t)\hat{\Phi}^{-1}(s)\hat{P}(s)=\hat{P}(t)\hat{\Phi}(t)\hat{\Phi}^{-1}(s),  \quad
\forall~ t,s\in I,\ee
and
\begin{eqnarray}
&&\E\|\hat{\Phi}(t)\hat{\Phi}^{-1}(s)\hat{P}(t)\|^2\leq \hat{M}
e^{-\hat{\alpha}(t-s)+\hat{\varepsilon} |s|}, \quad \forall~ (t,s)\in
I^{2}_{\geq},
\lb{c3}
\\
&&
\E\|\hat{\Phi}(t)\hat{\Phi}^{-1}(s)\hat{Q}(t)\|^2\leq \hat{M}
e^{-\hat{\alpha}(s-t)+\hat{\varepsilon} |s|}, \quad \forall~ (t,s)\in
I^{2}_{\leq},
\lb{x3}
\end{eqnarray}
where bound $\hat{M}:=40M$,  exponent $\hat{\alpha}:=\frac{\alpha}{2}-\frac{10M\tilde{M}}{\alpha}$, and nonuniform degree $\hat{\varepsilon}:=2\varepsilon$.
\end{theorem}

\mbox{\bf Proof of Theorem \ref{main23}.}  We first prove several lemmas which are essential in proving the theorem. The first one is the  existence and uniqueness lemma, which is slightly different from Lemma \ref{lem23} since $U(s,s)\xi_{0}$ is not necessarily equal to $\xi_{0}$ in \x{c4}. We will explain  the reason after Lemma \ref{lem38} under which condition there exists an equivalence between  \x{b9} and \x{c4} below.

\begin{lemma}\lb{lem32}
  For  any given initial value $\xi_{0} \in \R^{n}$, {\rm \x{a3}} has a unique solution $U(t,s) \xi_{0}$ with $U \in (\mathscr{L}_{c},
\|\cdot\|_{c})$ such that
   \bea\lb{c4} \EM U(t,s)=
\Phi(t)\Phi^{-1}(s)P(s)+ \int^{t}_{s}\Phi(t)\Phi^{-1}(\tau)P(\tau)H(\tau)U(\tau, s)d\omega(\tau)
\nn\\
\EM + \int^{t}_{s}\Phi(t)\Phi^{-1}(\tau)P(\tau)\tilde{B}(\tau)U(\tau, s)d\tau -\int^{\oo}_{t}\Phi(t)\Phi^{-1}(\tau)Q(\tau)H(\tau)U(\tau, s)d\omega(\tau)\nn\\
\EM -\int^{\oo}_{t}\Phi(t)\Phi^{-1}(\tau)Q(\tau)\tilde{B}(\tau)U(\tau, s)d\tau.
\eea
\end{lemma}
\prf{We first prove that the function $U(t,s) \xi_{0}$  is a solution of \x{a3}.
Set
\beaa \EM\xi(t)=\Phi^{-1}(s)P(s)\xi_{0}+ \int^{t}_{s}\Phi^{-1}(\tau)P(\tau)H(\tau)U(\tau, s)\xi_{0}d\omega(\tau)
\nn\\
\EM + \int^{t}_{s}\Phi^{-1}(\tau)P(\tau)\tilde{B}(\tau)U(\tau, s)\xi_{0}d\tau  -\int^{\oo}_{t}\Phi^{-1}(\tau)Q(\tau)H(\tau)U(\tau, s)\xi_{0}d\omega(\tau)\nn\\
\EM -\int^{\oo}_{t}\Phi^{-1}(\tau)Q(\tau)\tilde{B}(\tau)U(\tau, s)\xi_{0}d\tau.\eeaa
Let $y(t)=\Phi(t)\xi(t)$. Clearly,
\[ U(t,s) \xi_{0}= \Phi(t)\xi(t) =y(t),
\]
and then $\xi(t)$ satisfies the differential
\beaa d \xi(t)\EQ \Phi^{-1}(t)\big{(}B(t)-
G(t)H(t)\big{)}y(t)d t +\Phi^{-1}(t)H(t)y(t)d\omega(t).
\eeaa
Since $\Phi(t)$ is a fundamental matrix solution of \x{a1}. it follows from It\^{o} product rule that
\beaa dy(t)\EQ
d\Phi(t)\xi(t)+\Phi(t)d\xi(t)+G(t)\Phi(t)\Phi^{-1}(t)H(t)y(t)dt\\
\EQ A(t)y(t)dt+G(t)y(t)d\omega(t)+\big{(}B(t)-
G(t)H(t)\big{)}y(t)dt\\
\EM+H(t)y(t)d\omega(t)+G(t)H(t)y(t)dt\\\EQ
(A(t)+B(t))y(t)dt+(G(t)+H(t))y(t)d\omega(t),
\eeaa
which means that $y(t)=U(t,s)  \xi_{0}$ is a solution of \x{a3}.

Now we prove that $U$ is unique in $(\mathscr{L}_{c},
\|\cdot\|_{c})$.  Let
\beaa\lb{c5}\EM(\Gamma U)(t,s)=\Phi(t)\Phi^{-1}(s)P(s)+ \int^{t}_{s}\Phi(t)\Phi^{-1}(\tau)P(\tau)H(\tau)U(\tau, s)d\omega(\tau)
\nn\\
\EM + \int^{t}_{s}\Phi(t)\Phi^{-1}(\tau)P(\tau)\tilde{B}(\tau)U(\tau, s)d\tau -\int^{\oo}_{t}\Phi(t)\Phi^{-1}(\tau)Q(\tau)H(\tau)U(\tau, s)d\omega(\tau)\nn\\
\EM -\int^{\oo}_{t}\Phi(t)\Phi^{-1}(\tau)Q(\tau)\tilde{B}(\tau)U(\tau, s)d\tau.\eeaa
 The same idea as in Lemma \ref{lem23} can be applied to prove the uniqueness of the solution to \x{c4}. Squaring both sides of \x{c4}, and taking expectations, we have
\bea \EM \E\|(\Gamma U)(t,s)\|^{2} \le 5\E \|\Phi(t)\Phi^{-1}(s)P(s)\|^{2} +
5\E \left\|\int^{t}_{s}\Phi(t)\Phi^{-1}(\tau)P(\tau)H(\tau)U(\tau, s)d\omega(\tau)\right\|^{2}\nn\\
\EM + 5\E \left\|\int^{t}_{s}\Phi(t)\Phi^{-1}(\tau)Q(\tau)\tilde{B}(\tau)U(\tau, s)d\tau \right\|^{2}+ 5\E \left\|\int^{\oo}_{t}\Phi(t)\Phi^{-1}(\tau)Q(\tau)H(\tau)U(\tau, s)d\omega(\tau) \right\|^{2}\nn\\
\EM + 5\E \left\|\int^{\oo}_{t}\Phi(t)\Phi^{-1}(\tau)Q(\tau)\tilde{B}(\tau)U(\tau, s)d\tau \right\|^{2}\nn\\
\LE 5M e^{\varepsilon |s|} +10M e^{\varepsilon |s|} (\frac{\alpha h^{2}+8b^{2}+8g^{2}h^{2}}{\alpha^{2}})\sup_{(\tau,s)\in I^{2}_{\ge}}\left(\E\|U(\tau, s)\|^{2}e^{-\varepsilon |s|}\right),\nn
\eea
and this implies that
\[\E\|(\Gamma U)(t,s)\|^{2}e^{-\varepsilon s}\le 5M +\frac{10M\tilde{M}}{\alpha^{2}}\sup_{(\tau,s)\in I^{2}_{\ge}}\left(\E\|\hat{\Phi}(\tau, s)\|^{2}e^{-\varepsilon |s|}\right)< \oo\]
with $\tilde{M}=8b^{2}+8g^{2}h^{2}+\alpha h^2$.
Proceeding in the
same procedure as above, for any $U_{1}, U_{2}\in
\mathscr{L}_{c}$, we have
\be \lb{c6}\|\Gamma U_{1}-\Gamma U_{2}\|^{2}_{c}\le \frac{10M\tilde{M}}{\alpha^{2}}\sup_{(\tau,s)\in I^{2}_{\ge}}\left(\E\|U_{1}(\tau, s)-U_{2}(\tau, s)\|^{2}e^{-\varepsilon |s|}\right).\ee
Note that
\beaa\sup_{(t,s)\in I^{2}_{\ge}}\left(\E\|U_{1}(t, s)-U_{2}(t, s)\|^{2}e^{-\varepsilon s}\right)\LE
\sup_{(t,s)\in I^{2}_{\ge}}\left((\E\|U_{1}(t, s)-U_{2}(t, s)\|^{2})^{\frac{1}{2}}e^{-\frac{\varepsilon |s|}{2} }\right)^{2}\\
\LE \left(\sup_{(t,s)\in I^{2}_{\ge}}(\E\|U_{1}(t, s)-U_{2}(t, s)\|^{2})^{\frac{1}{2}}e^{-\frac{\varepsilon |s|}{2} }\right)^{2}\\
\EQ \|U_{1}-U_{2}\|^{2}_{c},
\eeaa
which together with \x{c6} implies
\[\|\Gamma U_{1}-\Gamma U_{2}\|_{c}\le \sqrt{\frac{10M\tilde{M}}{\alpha^{2}} } \|U_{1}-U_{2}\|_{c}.\]
Since $\tilde{M}<\frac{\alpha^{2}}{10M}$, $\Gamma$ is a contraction operator. Hence, there exists a
unique $U \in \mathscr{L}_{c}$ such that $\Gamma U= U$, which satisfies the identity \x{c4}. This completes the proof of the lemma. \hspace{\stretch{1}}$\Box$}

\begin{lemma}\lb{lem33}
  For  any $u\in (s,t)$ in $I$, we have
   \[U(t,s)=U(t,u)U(u,s)\]
in the sense of  $(\mathscr{L}_{c},
\|\cdot\|_{c})$.
\end{lemma}

\prf{By \x{c1} and \x{c4} with any $t\ge u \ge s$ in $I$, we have
  \bea\lb{c7} \EM U(t,u)U(u,s)=
\Phi(t)\Phi^{-1}(s)P(s)+ \int^{u}_{s}\Phi(t)\Phi^{-1}(\tau)P(\tau)H(\tau)U(\tau, s)d\omega(\tau)
\nn\\
\EM + \int^{u}_{s}\Phi(t)\Phi^{-1}(\tau)P(\tau)\tilde{B}(\tau)U(\tau, s)d\tau 
+\bigg{(}\int^{t}_{u}\Phi(t)\Phi^{-1}(\tau)P(\tau)H(\tau)U(\tau, u)d\omega(\tau)
\nn\\
\EM +\int^{t}_{u}\Phi(t)\Phi^{-1}(\tau)P(\tau)\tilde{B}(\tau)U(\tau, u)d\tau -\int^{\oo}_{t}\Phi(t)\Phi^{-1}(\tau)Q(\tau)H(\tau)U(\tau, u)d\omega(\tau)
\nn\\
\EM-\int^{\oo}_{t}\Phi(t)\Phi^{-1}(\tau)Q(\tau)\tilde{B}(\tau)U(\tau, u)d\tau\bigg{)}U(u,s).
\eea
Subtracting \x{c4} from \x{c7} we obtain
 \beaa\lb{c8} \EM U(t,s)-U(t,u)U(u,s)=
 \int^{t}_{u}\Phi(t)\Phi^{-1}(\tau)P(\tau)H(\tau)\left(U(\tau, s)-U(\tau,u)U(u,s)\right)d\omega(\tau)
\nn\\
\EM + \int^{t}_{u}\Phi(t)\Phi^{-1}(\tau)P(\tau)\tilde{B}(\tau)\left(U(\tau, s)-U(\tau,u)U(u,s)\right)d\tau\nn\\
\EM-\int^{\oo}_{t}\Phi(t)\Phi^{-1}(\tau)Q(\tau)H(\tau)\left(U(\tau, s)-U(\tau,u)U(u,s)\right)d\omega(\tau)
\nn\\
\EM-\int^{\oo}_{t}\Phi(t)\Phi^{-1}(\tau)Q(\tau)\tilde{B}(\tau)\left(U(\tau, s)-U(\tau,u)U(u,s)\right)d\tau.
\eeaa
Write $\tilde{U}(t,s)=U(t,s)-U(t,u)U(u,s)$. Now we prove $\tilde{U}$ is unique in $(\mathscr{L}_{c},
\|\cdot\|_{c})$. Let
\bea\lb{x12} \EM(\mathcal {T} \tilde{U})(t,s)= \int^{t}_{u}\Phi(t)\Phi^{-1}(\tau)P(\tau)H(\tau)\tilde{U}(\tau,s)d\omega(\tau)\nn\\
\EM +\int^{t}_{u}\Phi(t)\Phi^{-1}(\tau)P(\tau)\tilde{B}(\tau)\tilde{U}(\tau,s)d\tau
-\int^{\oo}_{t}\Phi(t)\Phi^{-1}(\tau)Q(\tau)H(\tau)\tilde{U}(\tau,s)d\omega(\tau)\nn\\
\EM-\int^{\oo}_{t}\Phi(t)\Phi^{-1}(\tau)Q(\tau)\tilde{B}(\tau)\tilde{U}(\tau,s).
\eea
Squaring both sides of \x{x12}, and taking expectations, it follows from \x{b10} that
\bea\lb{c9}  \E\|(\mathcal {T} \tilde{U})(t,s)\|^{2} \LE
4\E \left\|\int^{t}_{u}\Phi(t)\Phi^{-1}(\tau)P(\tau)H(\tau)\tilde{U}(\tau,s)d\omega(\tau)\right\|^{2}\nn\\
\EM + 4\E \left\| \int^{t}_{u}\Phi(t)\Phi^{-1}(\tau)P(\tau)\tilde{B}(\tau)\tilde{U}(\tau,s)d\tau \right\|^{2}\nn\\\EM+4\E \left\|\int^{\oo}_{t}\Phi(t)\Phi^{-1}(\tau)Q(\tau)H(\tau)\tilde{U}(\tau,s)d\omega(\tau)\right\|^{2}\nn\\
\EM + 4\E \left\| \int^{\oo}_{t}\Phi(t)\Phi^{-1}(\tau)Q(\tau)\tilde{B}(\tau)\tilde{U}(\tau,s)d\tau \right\|^{2}.
\eea
By using the It\^{o} isometry
property and the inequalities \x{c2}, the first term of the right-hand side in  \x{c9} can be deduced as follows:
\beaa \EM \E\left\|\int^{t}_{s}\Phi(t)\Phi^{-1}(\tau)P(\tau)H(\tau)\tilde{U}(\tau,s)
d\omega(\tau)\right\|^{2} \\
\EQ \int^{t}_{s}\E\|\Phi(t)\Phi^{-1}(\tau)P(\tau)\|^{2} \E\|H(\tau)\|^{2} \E\|\tilde{U}(\tau,s)\|^{2} d \tau\\
\LE Mh^{2}\int^{t}_{s}e^{-\alpha(t-\tau)}\E\|\tilde{U}(\tau,s)\|^{2} d \tau\\
\LE \frac{Mh^{2}}{\alpha}e^{\varepsilon |s|} \sup_{(\tau,s)\in I^{2}_{\ge}}\left(\E\|\tilde{U}(\tau, s)\|^{2}e^{-\varepsilon |s|}\right).
\eeaa

As to the second term in \x{c9}, it follows from
$\E\|x\|\le \sqrt{\E\|x\|^{2}}$, Cauchy-Schwarz inequality, It\^{o} isometry property of stochastic integrals, and \x{c2} that
\beaa \EM \E\left\|\int^{t}_{u}\Phi(t)\Phi^{-1}(\tau)P(\tau)\tilde{B}(\tau)\tilde{U}(\tau,s)d\tau \right\|^{2}\\
\LE \left(\int^{t}_{u}\E\left\|\Phi(t)\Phi^{-1}(\tau)P(\tau)\tilde{B}(\tau)\right\| d\tau \right)\\
\EM \times \left(\int^{t}_{u}\E\left\|\Phi(t)\Phi^{-1}(\tau)P(\tau)\tilde{B}(\tau)\right\| \E\left\|\tilde{U}(\tau,s)\right\|^{2}d\tau \right)\\
\LE 2M(b^{2}+g^{2}h^{2})\left(\int^{t}_{u}e^{-\frac{\alpha}{2}(t-\tau)} d\tau \right)
\left(\int^{t}_{u}e^{-\frac{\alpha}{2}(t-\tau)} \E\left\|\tilde{U}(\tau,s)\right\|^{2}d\tau \right)\\
\LE  \frac{8M(b^{2}+g^{2}h^{2})}{\alpha^{2}}e^{\varepsilon |s|}\sup_{(\tau,s)\in I^{2}_{\ge}}\left(\E\|\tilde{U}(\tau, s)\|^{2}e^{-\varepsilon |s|}\right).
\eeaa
Clearly, the proof above is also valid for proving the other terms in the right-hand
side in \x{c9}. Thus we can rewrite the inequality \x{c9} as
\beaa \EM \E\|(\mathcal {T} \tilde{U})(t,s)\|^{2} \le \frac{8M\tilde{M}}{\alpha^{2}}e^{\varepsilon |s|}\sup_{(\tau,s)\in I^{2}_{\ge}}\left(\E\|\tilde{U}(\tau, s)\|^{2}e^{-\varepsilon |s|}\right),
\eeaa
and
\beaa \|\mathcal {T} \tilde{U}\|_{c}\le \sqrt{\frac{8M\tilde{M}}{\alpha^{2}} }\|\tilde{U}\|_{c}
\eeaa
with $\tilde{M}=8b^{2}+8g^{2}h^{2}+\alpha h^2$. Proceeding in the
same procedure as above, for any $\tilde{U}_{1}, \tilde{U}_{2}\in
\mathscr{L}_{c}$, we have
\beaa \|\mathcal {T} \tilde{U}_{1}-\mathcal {T} \tilde{U}_{2}\|_{c}\le \sqrt{\frac{8M\tilde{M}}{\alpha^{2}} }\|\tilde{U}_{1}-\tilde{U}_{2}\|_{c}.
\eeaa
Since $\tilde{M}<\frac{\alpha^{2}}{8M}$, this implies $\mathcal {T}$ is a contraction. Hence, there is a unique $\tilde{U}\in (\mathscr{L}_{c},
\|\cdot\|_{c})$.
 On the other hand, $0\in (\mathscr{L}_{c},
\|\cdot\|_{c})$ also satisfies \x{x12}. Hence we must have
\[U(t,s)-U(t,u)U(u,s)=0\]
in $\mathscr{L}_{c}$.
Therefore, $U(t,s)=U(t,u)U(u,s)$ with $U \in (\mathscr{L}_{c},
\|\cdot\|_{c})$. This completes the proof of the lemma. \hspace{\stretch{1}}$\Box$}


\begin{lemma}\lb{lem34}Given $s\in I$, if $y(t):=\Lambda(t,s) \xi:[s,+\oo)\rightarrow L^{2}(\Omega, \R^{n})$ is a solution of {\rm\x{a3}} with $y(s)= \Lambda(s,s) \xi= \xi$ such that $\Lambda$ is bounded in  $(\mathscr{L}_{c},
\|\cdot\|_{c})$.
Then
   \bea\lb{c10} \EM y(t)=
\Phi(t)\Phi^{-1}(s)P(s)\xi+ \int^{t}_{s}\Phi(t)\Phi^{-1}(\tau)P(\tau)H(\tau)y(\tau)d\omega(\tau)
\nn\\
\EM + \int^{t}_{s}\Phi(t)\Phi^{-1}(\tau)P(\tau)\tilde{B}(\tau)y(\tau)d\tau -\int^{\oo}_{t}\Phi(t)\Phi^{-1}(\tau)Q(\tau)H(\tau)y(\tau)d\omega(\tau)\nn\\
\EM -\int^{\oo}_{t}\Phi(t)\Phi^{-1}(\tau)Q(\tau)\tilde{B}(\tau)y(\tau)d\tau.
\eea
\end{lemma}

\prf{It is easy to see from \x{b9} that
   \bea\lb{c11} P(t)y(t)\EQ
\Phi(t)\Phi^{-1}(s)P(s)\xi+ \int^{t}_{s}\Phi(t)\Phi^{-1}(\tau)P(\tau)H(\tau)y(\tau)d\omega(\tau)
\nn\\
\EM + \int^{t}_{s}\Phi(t)\Phi^{-1}(\tau)P(\tau)\tilde{B}(\tau)y(\tau)d\tau,\eea
and
\bea\lb{c12} Q(t)y(t)\EQ
\Phi(t)\Phi^{-1}(s)Q(s)\xi+ \int^{t}_{s}\Phi(t)\Phi^{-1}(\tau)Q(\tau)H(\tau)y(\tau)d\omega(\tau)
\nn\\
\EM + \int^{t}_{s}\Phi(t)\Phi^{-1}(\tau)Q(\tau)\tilde{B}(\tau)y(\tau)d\tau\eea
for each $(t,s)\in I^{2}_{\geq}$.
The equality \x{c12} can be rewritten in the equivalent form
\bea\lb{c13} Q(s)\xi\EQ
\Phi(s)\Phi^{-1}(t)Q(t)y(t)- \int^{t}_{s}\Phi(s)\Phi^{-1}(\tau)Q(\tau)H(\tau)y(\tau)d\omega(\tau)
\nn\\
\EM - \int^{t}_{s}\Phi(s)\Phi^{-1}(\tau)Q(\tau)\tilde{B}(\tau)y(\tau)d\tau.\eea
For convenience we can assume that $D=\|\Lambda\|_{c}<\oo$, since $\Lambda$ is bounded in  $(\mathscr{L}_{c},
\|\cdot\|_{c})$. Then it follows from \x{b8} and \x{x2} that
\beaa \E\|\Phi(s)\Phi^{-1}(t)Q(t)y(t)\|^{2}\le MD^{2}\|\xi\|^{2}e^{-\alpha(t-s)+\varepsilon(|t|+|s|)}.
\eeaa
Since $\alpha>\varepsilon$, the right hand side of this inequality goes to  zero as $t\rightarrow +\oo$. Furthermore, we have
\beaa \EM \E\left\|\int^{\oo}_{s}\Phi(s)\Phi^{-1}(\tau)Q(\tau)H(\tau)y(\tau)d\omega(\tau)\right\|^{2}\\
\EQ \int^{\oo}_{s} \E\left\|\Phi(s)\Phi^{-1}(\tau)Q(\tau)\right\|^{2}\E\left\|H(\tau)\right\|^{2}
\E\left\|y(\tau)\right\|^{2} d\tau\\
\LE \frac{h^{2}D^{2}M}{\alpha}e^{\varepsilon |s|},
\eeaa
and
\beaa \EM \E\left\|\int^{\oo}_{s}\Phi(s)\Phi^{-1}(\tau)Q(\tau)\tilde{B}(\tau)y(\tau)d\tau\right\|^{2}\\
\LE \left(\int^{\oo}_{s}\E\left\|\Phi(s)\Phi^{-1}(\tau)Q(\tau)\tilde{B}(\tau)\right\| d\tau \right) \left(\int^{\oo}_{s}\E\left\|\Phi(s)\Phi^{-1}(\tau)Q(\tau)\tilde{B}(\tau)\right\| \E\left\|y(\tau)\right\|^{2}d\tau \right)\\
\LE 2M(b^{2}+g^{2}h^{2})\left(\int^{\oo}_{s}e^{-\frac{\alpha}{2}(\tau-s)} d\tau \right)
\left(\int^{\oo}_{s}e^{-\frac{\alpha}{2}(\tau-s)} \E\left\|y(\tau)\right\|^{2}d\tau \right)\\
\LE \frac{8MD^{2}(b^{2}+g^{2}h^{2})}{\alpha^{2}}e^{\varepsilon |s|}.
\eeaa
Taking limits as $t\rightarrow +\oo$ in \x{c13}, we obtain
\beaa Q(s)\xi\EQ
- \int^{\oo}_{s}\Phi(s)\Phi^{-1}(\tau)Q(\tau)H(\tau)y(\tau)d\omega(\tau)
\nn\\
\EM - \int^{\oo}_{s}\Phi(s)\Phi^{-1}(\tau)Q(\tau)\tilde{B}(\tau)y(\tau)d\tau,\eeaa
and substitute it into \x{c12} yields
\beaa Q(t)y(t)\EQ
- \int^{\oo}_{s}\Phi(t)\Phi^{-1}(\tau)Q(\tau)H(\tau)y(\tau)d\omega(\tau)+ \int^{t}_{s}\Phi(t)\Phi^{-1}(\tau)Q(\tau)H(\tau)y(\tau)d\omega(\tau)
\nn\\
\EM - \int^{\oo}_{s}\Phi(t)\Phi^{-1}(\tau)Q(\tau)\tilde{B}(\tau)y(\tau)d\tau + \int^{t}_{s}\Phi(t)\Phi^{-1}(\tau)Q(\tau)\tilde{B}(\tau)y(\tau)d\tau\\
\EQ - \int^{\oo}_{t}\Phi(t)\Phi^{-1}(\tau)Q(\tau)H(\tau)y(\tau)d\omega(\tau)- \int^{\oo}_{t}\Phi(t)\Phi^{-1}(\tau)Q(\tau)\tilde{B}(\tau)y(\tau)d\tau.
\eeaa
Since $\xi$  is an arbitrary one in $\R^{n}$, then by adding this identity to \x{c11} yields the desired equation \x{c10}.  \hspace{\stretch{1}}$\Box$
}

Recall that $\hat{\Phi}(t,s)=\hat{\Phi}(t)\hat{\Phi}^{-1}(s)$ denotes the fundamental matrix solution of {\rm\x{a3}} with $\hat{\Phi}(s,s)=Id$. For each $t\in I$, define linear operators as
\be\lb{c14} \hat{P}(t)=\hat{\Phi}(t,t_{0})U(t_{0},t_{0})\hat{\Phi}(t_{0},t)\quad {\rm and}\quad
\hat{Q}(t)=Id-\hat{P}(t),
\ee
where $t_{0}$ is the left boundary point of the interval $I$. After presenting that
$\hat{P}(t)$ are projections, we prove the relationship \x{x4}, show the explicit expressions of the fundamental matrix solution $\hat{\Phi}(t,s)$ under the projections $\hat{P}(t)$, $\hat{Q}(t)$, and then deduce the inequalities \x{c3}-\x{x3}.
\begin{lemma}\lb{lem35}
  The operator $\hat{P}(t)$ are linear projections for $t\in I$, and {\rm \x{x4}} holds.
\end{lemma}
\prf{ By Lemma \ref{lem33}, we have $U(t_{0},t_{0})U(t_{0},t_{0})=U(t_{0},t_{0})$. Thus,
\[\hat{P}(t)\hat{P}(t)=\hat{\Phi}(t,t_{0})U(t_{0},t_{0})\hat{\Phi}(t_{0},t)
\hat{\Phi}(t,t_{0})U(t_{0},t_{0})\hat{\Phi}(t_{0},t)=\hat{P}(t).\]
Furthermore, for any $t, s\in I$, we obtain
\beaa \hat{P}(t)\hat{\Phi}(t,s)\EQ \hat{\Phi}(t,t_{0})U(t_{0},t_{0})\hat{\Phi}(t_{0},t)\hat{\Phi}(t,s)\\
\EQ \hat{\Phi}(t,s)\hat{\Phi}(s,t_{0})U(t_{0},t_{0})\hat{\Phi}(t_{0},s)=\hat{\Phi}(t,s)\hat{P}(s),
\eeaa
and this completes the proof of the lemma. \hspace{\stretch{1}}$\Box$
}

\begin{lemma}\lb{lem36}
 For  any given initial value $\xi_{0} \in \R^{n}$, the function  $\hat{P}(t)\hat{\Phi}(t,s)\xi_{0}$ is a solution of  {\rm \x{a3}} with $\hat{P}(t)\hat{\Phi}(t,s)$ is bounded  in  $(\mathscr{L}_{c},
\|\cdot\|_{c})$.
\end{lemma}
\prf{By Lemma \ref{lem32}, the function $U(t,t_{0})\xi_{0}$ is a solution of \x{a3} with initial value $U(t_{0},t_{0})\xi_{0}$ at time $t_{0}$. Clearly, $U(t,t_{0})=\hat{\Phi}(t,t_{0})U(t_{0},t_{0})$. Thus it is easy to see that
\beaa \hat{P}(t)\hat{\Phi}(t,s)\EQ \hat{\Phi}(t,t_{0})U(t_{0},t_{0})\hat{\Phi}(t_{0},t)\hat{\Phi}(t,s)
=U(t,t_{0})\hat{\Phi}(t_{0},s).
\eeaa
Therefore, it follows again from Lemma \ref{lem32} that $\hat{P}(t)\hat{\Phi}(t,s)\xi_{0}=U(t,t_{0})\hat{\Phi}(t_{0},s)\xi_{0}$ is a solution of \x{a3} with initial value $\hat{\Phi}(t_{0},s)\xi_{0} \in \R^{n}$. Moreover, from $U \in (\mathscr{L}_{c},
\|\cdot\|_{c})$ and the definition \x{b7}-\x{b8} of the space $(\mathscr{L}_{c},
\|\cdot\|_{c})$, we can see that $\hat{P}(t)\hat{\Phi}(t,s)$ is bounded in $(\mathscr{L}_{c},
\|\cdot\|_{c})$.
\hspace{\stretch{1}}$\Box$
}
\begin{lemma}\lb{lem38}
 For  any given initial value $\xi_{0} \in \R^{n}$, the function  $\hat{P}(t)\hat{\Phi}(t,s)\xi_{0}$ is a solution of  {\rm \x{a3}} with  $(t,s)\in I^{2}_{\geq}$ such that
\bea\lb{c15} \EM \hat{\Phi}(t,s)\hat{P}(s)=
\Phi(t)\Phi^{-1}(s)P(s)\hat{P}(s)+ \int^{t}_{s}\Phi(t)\Phi^{-1}(\tau)P(\tau)H(\tau)\hat{\Phi}(\tau, s)\hat{P}(s)d\omega(\tau)
\nn\\
\EM + \int^{t}_{s}\Phi(t)\Phi^{-1}(\tau)P(\tau)\tilde{B}(\tau)\hat{\Phi}(\tau, s)\hat{P}(s)d\tau -\int^{\oo}_{t}\Phi(t)\Phi^{-1}(\tau)Q(\tau)H(\tau)\hat{\Phi}(\tau, s)\hat{P}(s)d\omega(\tau)\nn\\
\EM -\int^{\oo}_{t}\Phi(t)\Phi^{-1}(\tau)Q(\tau)\tilde{B}(\tau)\hat{\Phi}(\tau, s)\hat{P}(s)d\tau.
\eea
\end{lemma}
\prf{Let $y(t)=\hat{P}(t)\hat{\Phi}(t,s)\xi_{0}$ with given $s\in I$, and denote $\xi=\hat{P}(s)\xi_{0}$ the initial condition at time $s$.  Clearly, $y(t)$ is a solution of \x{a3} with $y(s)=\hat{P}(s)\xi_{0}=\hat{P}(s)\hat{P}(s)\xi_{0}=\xi$. By Lemma \ref{lem36}, $\hat{P}(t)\hat{\Phi}(t,s)$ is bounded  in  $(\mathscr{L}_{c},
\|\cdot\|_{c})$. Since $\xi_{0}$ is arbitrary in $\R^{n}$, the identity \x{c15} follows now readily from Lemma \ref{lem34}. \hspace{\stretch{1}}$\Box$
}

\begin{remark}
  From Lemma {\rm \ref{lem38}}, we know that the explicit expressions {\rm\x{b9}} and {\rm\x{c4}} are the same under the condition of NMS-EC. In fact,  as a special case of Lemma {\rm\ref{lem38}}, $\hat{\Phi}(t,s)$ is always  bounded  in  $(\mathscr{L}_{c},
\|\cdot\|_{c})$ with $I=[t_{0},+\oo)$ since projections are the identity.
\end{remark}

 In the following lemma, we present the explicit expression of $\hat{\Phi}(t,s)\hat{Q}(s)$ with  $(t,s)\in I^{2}_{\leq}$.

\begin{lemma}\lb{lem37}
 For  any given initial value $\xi_{0} \in \R^{n}$, the function  $\hat{Q}(t)\hat{\Phi}(t,s)\xi_{0}$ is a solution of  {\rm \x{a3}} with  $(t,s)\in I^{2}_{\leq}$ such that
\bea\lb{c18} \EM \hat{\Phi}(t,s)\hat{Q}(s)=
\Phi(t)\Phi^{-1}(s)Q(s)\hat{Q}(s)+ \int^{t}_{t_{0}}\Phi(t)\Phi^{-1}(\tau)P(\tau)H(\tau)\hat{\Phi}(\tau, s)\hat{Q}(s)d\omega(\tau)
\nn\\
\EM + \int^{t}_{t_{0}}\Phi(t)\Phi^{-1}(\tau)P(\tau)\tilde{B}(\tau)\hat{\Phi}(\tau, s)\hat{Q}(s)d\tau  -\int^{s}_{t}\Phi(t)\Phi^{-1}(\tau)Q(\tau)H(\tau)\hat{\Phi}(\tau, s)\hat{Q}(s)d\omega(\tau)\nn\\
\EM -\int^{s}_{t}\Phi(t)\Phi^{-1}(\tau)Q(\tau)\tilde{B}(\tau)\hat{\Phi}(\tau, s)\hat{Q}(s)d\tau.
\eea
\end{lemma}
\prf{Following the same lines as given in the proof of Lemma \ref{lem23}, one can prove that
   \beaa\lb{c19} \hat{\Phi}(t,s) \EQ
\Phi(t)\Phi^{-1}(s)+ \int^{t}_{s}\Phi(t)\Phi^{-1}(\tau)H(\tau)\hat{\Phi}(\tau, s)d\omega(\tau)
\nn\\
\EM + \int^{t}_{s}\Phi(t)\Phi^{-1}(\tau)\tilde{B}(\tau)\hat{\Phi}(\tau, s)d\tau\eeaa
for any $(t,s) \in I^{2}_{\leq}$. Write $K(t)=\hat{\Phi}(t,t_{0})\hat{Q}(t_{0})$.
Therefore,  \bea\lb{c20} K(t)\EQ
\Phi(t)\Phi^{-1}(t_{0})\hat{Q}(t_{0})+ \int^{t}_{t_{0}}\Phi(t)\Phi^{-1}(\tau)H(\tau)K(\tau)d\omega(\tau)
\nn\\
\EM + \int^{t}_{t_{0}}\Phi(t)\Phi^{-1}(\tau)\tilde{B}(\tau)K(\tau)d\tau.\eea
On the other hand, it follows from $\hat{P}(t)=\hat{\Phi}(t,t_{0})U(t_{0},t_{0})\hat{\Phi}(t_{0},t)$ and \x{c4} with $t=s=t_{0}$ that
 \bea\lb{c21} \hat{P}(t_{0}) \EQ U(t_{0},t_{0})=
P(t_{0})
 -\int^{\oo}_{t_{0}}\Phi(t_{0})\Phi^{-1}(\tau)Q(\tau)H(\tau)U(\tau, t_{0})d\omega(\tau)\nn\\
\EM -\int^{\oo}_{t_{0}}\Phi(t_{0})\Phi^{-1}(\tau)Q(\tau)\tilde{B}(\tau)U(\tau, t_{0})d\tau.
\eea
Since $P(t_{0})$ and  $Q(t_{0})$ are complementary projections, multiplies \x{c21} on the left with $P(t_{0})$. This gives
\bea\lb{c22} P(t_{0})\hat{P}(t_{0})=P(t_{0}).
\eea
In addition,
\be\lb{c23} Q(t_{0})\hat{Q}(t_{0})=\big(Id-P(t_{0})\big)\big(Id-\hat{P}(t_{0})\big)
=Id-\hat{P}(t_{0})=\hat{Q}(t_{0}).
\ee
By \x{c20}, using \x{c23}, we have
\beaa \Phi(t)\Phi^{-1}(s)Q(s)K(s)\EQ
\Phi(t)\Phi^{-1}(t_{0})\hat{Q}(t_{0})+ \int^{s}_{t_{0}}\Phi(t)\Phi^{-1}(\tau)Q(\tau)H(\tau)K(\tau)d\omega(\tau)
\nn\\
\EM + \int^{s}_{t_{0}}\Phi(t)\Phi^{-1}(\tau)Q(\tau)\tilde{B}(\tau)K(\tau)d\tau,\eeaa
which can be rewritten as
\bea\lb{c24} \Phi(t)\Phi^{-1}(t_{0})\hat{Q}(t_{0})\EQ \Phi(t)\Phi^{-1}(s)Q(s)K(s)
- \int^{s}_{t_{0}}\Phi(t)\Phi^{-1}(\tau)Q(\tau)H(\tau)K(\tau)d\omega(\tau)
\nn\\
\EM - \int^{s}_{t_{0}}\Phi(t)\Phi^{-1}(\tau)Q(\tau)\tilde{B}(\tau)K(\tau)d\tau.\eea
Substitute \x{c24} into \x{c20} leads to
\bea\lb{c25} \EM K(t)=
\Phi(t)\Phi^{-1}(s)Q(s)K(s)
- \int^{s}_{t_{0}}\Phi(t)\Phi^{-1}(\tau)Q(\tau)H(\tau)K(\tau)d\omega(\tau)
\nn\\
\EM - \int^{s}_{t_{0}}\Phi(t)\Phi^{-1}(\tau)Q(\tau)\tilde{B}(\tau)K(\tau)d\tau+ \int^{t}_{t_{0}}\Phi(t)\Phi^{-1}(\tau)H(\tau)K(\tau)d\omega(\tau)
\nn\\
\EM + \int^{t}_{t_{0}}\Phi(t)\Phi^{-1}(\tau)\tilde{B}(\tau)K(\tau)d\tau\nn \\
\EM = \Phi(t)\Phi^{-1}(s)Q(s)K(s)
- \int^{s}_{t}\Phi(t)\Phi^{-1}(\tau)Q(\tau)H(\tau)K(\tau)d\omega(\tau)
\nn\\
\EM - \int^{s}_{t}\Phi(t)\Phi^{-1}(\tau)Q(\tau)\tilde{B}(\tau)K(\tau)d\tau+ \int^{t}_{t_{0}}\Phi(t)\Phi^{-1}(\tau)P(\tau)H(\tau)K(\tau)d\omega(\tau)
\nn\\
\EM + \int^{t}_{t_{0}}\Phi(t)\Phi^{-1}(\tau)P(\tau)\tilde{B}(\tau)K(\tau)d\tau.\eea
Since \x{x4} we have $K(t)=\hat{\Phi}(t,t_{0})\hat{Q}(t_{0})=\hat{Q}(t)\hat{\Phi}(t,t_{0})$. Therefore,
$K(t)\hat{\Phi}(t_{0},s)=\hat{Q}(t)\hat{\Phi}(t,s)$ for every $(t,s) \in I^{2}_{\le}$. Thus, multiplying \x{c25} on the right with $\hat{\Phi}(t_{0},s)$. This yields the desired identity
\x{c18}. \hspace{\stretch{1}}$\Box$
}

We proceed with the proof of Theorem \ref{main23}.
 Squaring both sides of  \x{c15}, and taking expectations. Setting $z(t,s)=\E\|\hat{\Phi}(t,s)\hat{P}(s)\|^{2}$ with $(t,s) \in I^{2}_{\ge}$.
  It follows from \x{b10} that
\bea\lb{c26}  z(t,s) \LE
5\E\|\Phi(t)\Phi^{-1}(s)P(s)\hat{P}(s)\|^{2}+ 5\E\left\|\int^{t}_{s}\Phi(t)\Phi^{-1}(\tau)P(\tau)H(\tau)\hat{\Phi}(\tau, s)\hat{P}(s)d\omega(\tau)\right\|^{2} \nn\\
\EM+ 5\E\left\|\int^{t}_{s}\Phi(t)\Phi^{-1}(\tau)P(\tau)\tilde{B}(\tau)\hat{\Phi}(\tau, s)\hat{P}(s)d\tau\right\|^{2}\nn\\
\EM+ 5\E\left\|\int^{\oo}_{t}\Phi(t)\Phi^{-1}(\tau)Q(\tau)H(\tau)\hat{\Phi}(\tau, s)\hat{P}(s)d\omega(\tau)\right\|^{2}\nn\\
\EM+ 5\E\left\|\int^{\oo}_{t}\Phi(t)\Phi^{-1}(\tau)Q(\tau)\tilde{B}(\tau)\hat{\Phi}(\tau, s)\hat{P}(s)d\tau\right\|^{2}.
\eea
By using the It\^{o} isometry
property and the  inequalities \x{c2}, the second term of right-hand side in  \x{c26} can be deduced as follows:
\beaa \EM \E\left\|\int^{t}_{s}\Phi(t)\Phi^{-1}(\tau)P(\tau)H(\tau)\hat{\Phi}(\tau, s)\hat{P}(s)d\omega(\tau)\right\|^{2} \\
\EQ \int^{t}_{s}\E\|\Phi(t)\Phi^{-1}(\tau)P(\tau)\|^{2} \E\|H(\tau)\|^{2} \E\|\hat{\Phi}(\tau, s)\hat{P}(s)\|^{2} d \tau\\
\LE Mh^{2}\int^{t}_{s}e^{-\alpha(t-\tau)}\E\|\hat{\Phi}(\tau,s)\hat{P}(s)\|^{2} d \tau.
\eeaa
As to the third term in \x{c26}, it follows from $\E\|x\|\le \sqrt{\E\|x\|^{2}}$, Cauchy-Schwarz inequality, and the inequalities \x{c2} that
\beaa \EM \E\left\|\int^{t}_{s}\Phi(t)\Phi^{-1}(\tau)P(\tau)\tilde{B}(\tau)\hat{\Phi}(\tau, s)\hat{P}(s)d\tau\right\|^{2}\\
\EQ \E\left\|\int^{t}_{s}\left(\Phi(t)\Phi^{-1}(\tau)P(\tau)\tilde{B}(\tau)\right)^{\frac{1}{2}}
\left(\left(\Phi(t)\Phi^{-1}(\tau)P(\tau)\tilde{B}(\tau)\right)^{\frac{1}{2}}\hat{\Phi}(\tau,s)\hat{P}(s)\right)
d\tau\right\|^{2}\\
\LE \left(\int^{t}_{s}\E\left\|\Phi(t)\Phi^{-1}(\tau)P(\tau)\tilde{B}(\tau)\right\| d\tau \right)\\
\EM \times \left(\int^{t}_{s}\E\left\|\Phi(t)\Phi^{-1}(\tau)P(\tau)\tilde{B}(\tau)\right\| \E\left\|\hat{\Phi}(\tau,s)\hat{P}(s)\right\|^{2}d\tau \right)\\
\LE 2M(b^{2}+g^{2}h^{2})\left(\int^{t}_{s}e^{-\frac{\alpha}{2}(t-\tau)} d\tau \right)
\left(\int^{t}_{s}e^{-\frac{\alpha}{2}(t-\tau)} \E\left\|\hat{\Phi}(\tau,s)\hat{P}(s)\right\|^{2}d\tau \right)\\
\LE  \frac{4M(b^{2}+g^{2}h^{2})}{\alpha}\int^{t}_{s}e^{-\frac{\alpha}{2}(t-\tau)} \E\left\|\hat{\Phi}(\tau,s)\hat{P}(s)\right\|^{2}d\tau.
\eeaa
Clearly, the proof above is also valid for proving the other terms in the right-hand
side in \x{c26}. Thus we can rewrite the inequality \x{c26} as
\bea\lb{c27} z(t,s)\LE 5Me^{-\alpha(t-s)+\varepsilon |s|}z(s,s)+5Mh^{2}\left(\int^{t}_{s}e^{-\alpha(t-\tau)}z(\tau,s) d \tau+
\int^{\oo}_{t}e^{-\alpha(\tau-t)}z(\tau,s) d \tau\right)\nn\\
\EM +\frac{20M(b^{2}+g^{2}h^{2})}{\alpha}\left(\int^{t}_{s}e^{-\frac{\alpha}{2}(t-\tau)} z(\tau,s)d\tau
+\int^{\oo}_{t}e^{-\frac{\alpha}{2}(\tau-t)} z(\tau,s)d\tau
\right)\nn\\
\LE 5Me^{-\frac{\alpha}{2}(t-s)+\varepsilon |s|}z(s,s)+ \frac{5M\tilde{M}}{\alpha}\left(\int^{t}_{s}e^{-\frac{\alpha}{2}(t-\tau)}z(\tau,s) d \tau+
\int^{\oo}_{t}e^{-\frac{\alpha}{2}(\tau-t)}z(\tau,s) d \tau\right)\nn\\
\eea
with $\tilde{M}=8b^{2}+8g^{2}h^{2}+\alpha h^2$. Let
\[ Z(t,s)=5Me^{-\frac{\alpha}{2}(t-s)+\varepsilon |s|}z(s,s)+ \frac{5M\tilde{M}}{\alpha}\left(\int^{t}_{s}e^{-\frac{\alpha}{2}(t-\tau)}z(\tau,s) d \tau+
\int^{\oo}_{t}e^{-\frac{\alpha}{2}(\tau-t)}z(\tau,s) d \tau\right).
\]
 Clearly, inequality \x{c27} can be rewritten as
\[z(t,s)\le Z(t,s).\]
On the other hand,
\[\frac{d}{dt}Z(t,s)=-\frac{\alpha}{2}Z(t,s)+\frac{10M \tilde{M}}{\alpha}z(t,s),\] and therefore, \[\frac{d}{dt}Z(t,s)\leq \left(\frac{10M \tilde{M}}{\alpha}-\frac{\alpha}{2}\right) Z(t,s).\]
Integrating the above inequality from $s$ to $t$ and note that
$Z(s,s)=5Me^{\varepsilon |s|}z(s,s)$, we obtain \[z(t,s)\leq Z(t,s)\leq 5Me^{\varepsilon |s|}
e^{(-\frac{\alpha}{2}+\frac{10M\tilde{M}}{\alpha})(t-s)}z(s,s),\quad \forall~ (t,s)\in
I^{2}_{\geq}.\]
By $z(t,s)=\E\|\hat{\Phi}(t,s)\hat{P}(s)\|^{2}$, we have
\be\lb{c28}\E\|\hat{\Phi}(t,s)\hat{P}(s)\|^{2}\leq 5M
e^{(-\frac{\alpha}{2}+\frac{10M\tilde{M}}{\alpha})(t-s)+\varepsilon |s|}\E\|\hat{P}(s)\|^{2},\quad \forall~ (t,s)\in
I^{2}_{\geq}.\ee

Similarly, squaring both sides of  \x{c18}, and taking expectations. Using the same way as above, we obtain
\be\lb{c29}\E\|\hat{\Phi}(t,s)\hat{Q}(s)\|^{2}\leq 5M
e^{(-\frac{\alpha}{2}+\frac{10M\tilde{M}}{\alpha})(s-t)+\varepsilon |s|}\E\|\hat{Q}(s)\|^{2},\quad \forall~ (t,s)\in
I^{2}_{\leq}.\ee

Now we try to find out the bounds in mean square setting for the projections $\hat{P}(t)$,  $\hat{Q}(t)$. Multiplying \x{c15} with $Q(t)$ on the left side, and let $t=s$, we have
\bea\lb{c30}  Q(t)\hat{P}(t)\EQ
 -\int^{\oo}_{t}\Phi(t)\Phi^{-1}(\tau)Q(\tau)H(\tau)\hat{\Phi}(\tau, t)\hat{P}(t)d\omega(\tau)\nn\\
\EM -\int^{\oo}_{t}\Phi(t)\Phi^{-1}(\tau)Q(\tau)\tilde{B}(\tau)\hat{\Phi}(\tau, t)\hat{P}(t)d\tau.
\eea
 By \x{c30}, using \x{x2}, \x{x1} and \x{c28},  we have
\bea\lb{c31}\E\|Q(t)\hat{P}(t)\|^{2}\LE 2\E \left\|\int^{\oo}_{t}\Phi(t)\Phi^{-1}(\tau)Q(\tau)H(\tau)\hat{\Phi}(\tau, t)\hat{P}(t)d\omega(\tau) \right\|^{2}\nn\\
\EM +2\E \left\|\int^{\oo}_{t}\Phi(t)\Phi^{-1}(\tau)Q(\tau)\tilde{B}(\tau)\hat{\Phi}(\tau, t)\hat{P}(t)d\tau \right\|^{2}\nn\\
\LE 2\int^{\oo}_{t} \E\|\Phi(t)\Phi^{-1}(\tau)Q(\tau)\|^{2}
\E\|H(\tau)\|^{2} \E\|\hat{\Phi}(\tau, t)\hat{P}(t)\|^{2} d\tau \nn\\
\EM + 2 \left(\int^{\oo}_{t} \E \|\Phi(t)\Phi^{-1}(\tau)Q(\tau)\|\E\|\tilde{B}(\tau)\|^{\frac{1}{2}}d \tau\right)\nn\\
 \EM\times \left(\int^{\oo}_{t} \E \|\Phi(t)\Phi^{-1}(\tau)Q(\tau)\|\E\|\tilde{B}(\tau)\|^{\frac{3}{2}} \E\|\hat{\Phi}(\tau, t)\hat{P}(t)\|^{2}d \tau\right)\nn\\
 \LE \frac{10M^{2}\tilde{M}}{\alpha}\E\|\hat{P}(t)\|^{2}\int^{\oo}_{t}
 e^{-(\alpha+\tilde{\alpha}-\varepsilon)(\tau-t)}d \tau\nn\\
 \LE \frac{10M^{2}\tilde{M}}{\alpha(\alpha+\tilde{\alpha}-\varepsilon)}\E\|\hat{P}(t)\|^{2},
\eea
since $\alpha>\varepsilon$ and $\tilde{\alpha}=\frac{\alpha}{2}-\frac{10M\tilde{M}}{\alpha}>0$.
In addition, it follows from \x{c18} with $t=s$ that
\bea\lb{c33} \EM P(t)\hat{Q}(t)=
\int^{t}_{t_{0}}\Phi(t)\Phi^{-1}(\tau)P(\tau)H(\tau)\hat{\Phi}(\tau, t)\hat{Q}(t)d\omega(\tau)
\nn\\
\EM + \int^{t}_{t_{0}}\Phi(t)\Phi^{-1}(\tau)P(\tau)\tilde{B}(\tau)\hat{\Phi}(\tau, t)\hat{Q}(t)d\tau.
\eea
Similarly, by \x{c33}, using \x{c2}, \x{x1} and \x{c29}, we obtain
\be\lb{c34} \E\|P(t)\hat{Q}(t)\|^{2}\le \frac{10M^{2}\tilde{M}}{\alpha(\alpha+\tilde{\alpha}-\varepsilon)}\E\|\hat{Q}(t)\|^{2}.
\ee
Meanwhile, notice that
\beaa \hat{P}(t)-P(t) \EQ \hat{P}(t)-P(t)\hat{P}(t)-P(t)+P(t)\hat{P}(t)\nn\\
\EQ (Id-P(t))\hat{P}(t)-(Id-\hat{P}(t))P(t)\nn\\
\EQ Q(t)\hat{P}(t)-P(t)\hat{Q}(t).
\eeaa
Thus it follows from \x{c31} and \x{c34} that
\be\lb{c42}\E\|\hat{P}(t)-P(t)\|^{2}\le \frac{20M^{2}\tilde{M}}{\alpha(\alpha+\tilde{\alpha}-\varepsilon)}
(\E\|\hat{P}(t)\|^{2}+\E\|\hat{Q}(t)\|^{2}).
\ee
On the other hand, it follows from \x{c2}-\x{x2} with $t=s$ that
\[ \E\|P(t)\|^{2}\le M e^{\varepsilon |t|},\quad {\rm and}\quad \E\|Q(t)\|^{2}\le M e^{\varepsilon |t|}.\]
Therefore,
\beaa \E\|\hat{P}(t)\|^{2}\LE 2\E\|\hat{P}(t)-P(t)\|^{2}+2\E\|P(t)\|^{2}\\
\LE
\frac{40M^{2}\tilde{M}}{\alpha(\alpha+\tilde{\alpha}-\varepsilon)}
(\E\|\hat{P}(t)\|^{2}+\E\|\hat{Q}(t)\|^{2})+2M e^{\varepsilon |t|}.\eeaa
Since $\hat{Q}(t)-Q(t)=(Id-\hat{P}(t))-(Id-P(t))=P(t)-\hat{P}(t)$, we also have
\beaa \E\|\hat{Q}(t)\|^{2}\LE 2\E\|\hat{P}(t)-P(t)\|^{2}+2\E\|Q(t)\|^{2}\\
\LE
\frac{40M^{2}\tilde{M}}{\alpha(\alpha+\tilde{\alpha}-\varepsilon)}
(\E\|\hat{P}(t)\|^{2}+\E\|\hat{Q}(t)\|^{2})+2M e^{\varepsilon |t|}.\eeaa
Then we know
\[(\E\|\hat{P}(t)\|^{2}+\E\|\hat{Q}(t)\|^{2})\le \frac{80M^{2}\tilde{M}}{\alpha(\alpha+\tilde{\alpha}-\varepsilon)}
(\E\|\hat{P}(t)\|^{2}+\E\|\hat{Q}(t)\|^{2})+4M e^{\varepsilon |t|},
\]
and hence,
\[\left(1-\frac{80M^{2}\tilde{M}}{\alpha(\alpha+\tilde{\alpha}-\varepsilon)}\right)
(\E\|\hat{P}(t)\|^{2}+\E\|\hat{Q}(t)\|^{2})\le
4M e^{\varepsilon |t|}.
\]
Since $\tilde{M}:=8b^{2}+8g^{2}h^{2}+\alpha h^2$, we can obtain
\[\frac{80M^{2}\tilde{M}}{\alpha(\alpha+\tilde{\alpha}-\varepsilon)}\le \frac{1}{2}\] by letting $b$ and $h$ sufficiently small. This yields
\be\lb{c35} \E\|\hat{P}(t)\|^{2}\le 8M e^{\varepsilon |t|}\quad {\rm and}\quad \E\|\hat{Q}(t)\|^{2}\le 8M e^{\varepsilon |t|}.\ee
By \x{c28}, \x{c29}, using \x{c35} we obtain
\[\E\|\hat{\Phi}(t,s)\hat{P}(s)\|^{2}\leq 40M^{2}
e^{(-\frac{\alpha}{2}+\frac{10M\tilde{M}}{\alpha})(t-s)+2\varepsilon |s|},\quad \forall~ (t,s)\in
I^{2}_{\geq},\]
and
\[\E\|\hat{\Phi}(t,s)\hat{Q}(s)\|^{2}\leq 40M
e^{(-\frac{\alpha}{2}+\frac{10M\tilde{M}}{\alpha})(s-t)+2\varepsilon |s|},\quad \forall~ (t,s)\in
I^{2}_{\leq}.\]
This completes the proof of the theorem. \hspace{\stretch{1}}$\Box$

Under the hypotheses of Theorem {\rm\ref{main23}}, the following theorem try to discuss the differences of projections $P(t)$ and
$\hat{P}(t)$ in the mean square sense.
To illustrate it clearly, write \[\Phi(t,s)=\Phi(t)\Phi^{-1}(s).\] Obviously, $\Phi(t,s)$  is a fundamental matrix solution of {\rm\x{a1}}  with $\Phi(s,s)=Id$.

\begin{theorem} \lb{main32}
  Under the hypotheses of Theorem {\rm\ref{main23}}, for any $t\in I$, we have
  \be\lb{c39} P(t)=\Phi(t_{0},t)P(t_{0})\Phi(t,t_{0}), \quad {\rm and} \quad \hat{P}(t)=\hat{\Phi}(t_{0},t)\hat{P}(t_{0})\hat{\Phi}(t,t_{0}),
  \ee
  and
  \be\lb{c40}\E\|P(t)-\hat{P}(t)\|^{2}\le \frac{320M^{3}\tilde{M}}{\alpha(\alpha+\tilde{\alpha}-\varepsilon)}e^{\varepsilon |t|}.
  \ee
  In particular, for each fixed $t\in I$, we have $\E\|P(t)-\hat{P}(t)\|^{2}\rightarrow 0$ as $b, h \rightarrow 0$.
\end{theorem}

\prf{The second equality of \x{c39} is obvious from the definition \x{c14} of linear operators $\hat{P}(t)$. For the first term in \x{c39}, it follows from \x{c1} that
\[P(t)\Phi(t,t_{0})\Phi(t_{0},s)=\Phi(t,t_{0})\Phi(t_{0},s)P(s),  \quad
\forall~ t,s\in I,\]
and then
\be\lb{c41}\Phi(t_{0},t)P(t)\Phi(t,t_{0})=\Phi(t_{0},s)P(s)\Phi(s,t_{0}),  \quad
\forall~ t,s\in I.\ee
Taking $s=t_{0}$ in \x{c41}, we obtain
\[\Phi(t_{0},t)P(t)\Phi(t,t_{0})=P(t_{0}).\]
Thus,
\[P(t)=\Phi(t_{0},t)P(t_{0})\Phi(t,t_{0}).\]
In addition, \x{c40} follows immediately from \x{c42} and \x{c35}. \hspace{\stretch{1}}$\Box$
}

\begin{theorem}\lb{main33}
   Under the hypotheses of Theorem {\rm\ref{main23}}, we have
   \[ \E\|\hat{\Phi}(t,s)\hat{P}(s)-\Phi(t,s)P(s)\|^{2}\le \frac{720M\tilde{M}}{\alpha-\hat{\alpha}}e^{-\hat{\alpha}(t-s)+\hat{\varepsilon}|s|},
   \quad\forall~ (t,s)\in
I^{2}_{\geq},
   \]
   and
   \[ \E\|\hat{\Phi}(t,s)\hat{Q}(s)-\Phi(t,s)Q(s)\|^{2}\le \frac{720M\tilde{M}}{\alpha-\hat{\alpha}}e^{-\hat{\alpha}(s-t)+\hat{\varepsilon}|s|},
   \quad\forall~ (t,s)\in
I^{2}_{\leq}.
   \]
\end{theorem}

\prf{By $\hat{P}(s)\hat{P}(s)=\hat{P}(s)$, it follows from \x{c15} that
\bea\lb{c44}
\EM
\E \left\|\hat{\Phi}(t,s)\hat{P}(s)\hat{P}(s)-\Phi(t,s)P(s)\hat{P}(s)\right\|^{2}
\le 4\E \left\|\int^{t}_{s}\Phi(t,\tau)P(\tau)H(\tau)\hat{\Phi}(\tau, s)\hat{P}(s)d\omega(\tau)\right\|^{2}\nn\\
\EM +4\E \left\|\int^{t}_{s}\Phi(t,\tau)P(\tau)\tilde{B}(\tau)\hat{\Phi}(\tau, s)\hat{P}(s)d\tau \right\|^{2} +4\E \left\|\int^{\oo}_{t}\Phi(t,\tau)Q(\tau)H(\tau)\hat{\Phi}(\tau, s)\hat{P}(s)d\omega(\tau) \right\|^{2}\nn\\
\EM +4\E \left\|\int^{\oo}_{t}\Phi(t,\tau)Q(\tau)\tilde{B}(\tau)\hat{\Phi}(\tau, s)\hat{P}(s)d\tau\right\|^{2}.
\eea

By \x{c2} and \x{c3}, using $\alpha-\hat{\alpha}=\frac{\alpha}{2}+\frac{10M\tilde{M}}{\alpha}>0$, the first term of right-hand side in  \x{c44} can be deduced as follows:
\beaa \EM \E\left\|\int^{t}_{s}\Phi(t)\Phi^{-1}(\tau)P(\tau)H(\tau)\hat{\Phi}(\tau, s)\hat{P}(s)d\omega(\tau)\right\|^{2} \\
\EQ \int^{t}_{s}\E\|\Phi(t)\Phi^{-1}(\tau)P(\tau)\|^{2} \E\|H(\tau)\|^{2} \E\|\hat{\Phi}(\tau, s)\hat{P}(s)\|^{2} d \tau\\
\LE M\hat{M}h^{2}\int^{t}_{s}e^{-\alpha(t-\tau)}
e^{-\hat{\alpha}(\tau-s)+\hat{\varepsilon} |s|} d \tau\\
\EQ M\hat{M}h^{2}e^{-\alpha(t-s)+\hat{\varepsilon} |s|}\int^{t}_{s}e^{(\alpha-\hat{\alpha})(\tau-s)}d \tau\\
\LE \frac{M\hat{M}h^{2}}{\alpha-\hat{\alpha}}e^{-\hat{\alpha}(t-s)+\hat{\varepsilon} |s|}.
\eeaa
As to the second term in \x{c44}, by $\frac{\alpha}{2}-\hat{\alpha}=\frac{10M\tilde{M}}{\alpha}>0$, we have $2\alpha^{2}-\alpha\hat{\alpha}>0$. It follows from $\E\|x\|\le \sqrt{\E\|x\|^{2}}$, Cauchy-Schwarz inequality, and the inequalities \x{c2}, \x{c3} that
\beaa \EM \E\left\|\int^{t}_{s}\Phi(t,\tau)P(\tau)\tilde{B}(\tau)\hat{\Phi}(\tau, s)\hat{P}(s)d\tau\right\|^{2}\\
\EQ \E\left\|\int^{t}_{s}\left(\Phi(t,\tau)P(\tau)\tilde{B}(\tau)\right)^{\frac{1}{2}}
\left(\left(\Phi(t,\tau)P(\tau)\tilde{B}(\tau)\right)^{\frac{1}{2}}\hat{\Phi}(\tau,s)\hat{P}(s)\right)
d\tau\right\|^{2}\\
\LE \left(\int^{t}_{s}\E\left\|\Phi(t,\tau)P(\tau)\tilde{B}(\tau)\right\| d\tau \right)\\
\EM \times \left(\int^{t}_{s}\E\left\|\Phi(t,\tau)P(\tau)\tilde{B}(\tau)\right\| \E\left\|\hat{\Phi}(\tau,s)\hat{P}(s)\right\|^{2}d\tau \right)\\
\LE 2M\hat{M}(b^{2}+g^{2}h^{2})\left(\int^{t}_{s}e^{-\frac{\alpha}{2}(t-\tau)} d\tau \right)
\left(\int^{t}_{s}e^{-\frac{\alpha}{2}(t-\tau)} e^{-\hat{\alpha}(\tau-s)+\hat{\varepsilon} |s|}d\tau \right)\\
\LE  \frac{8M\hat{M}(b^{2}+g^{2}h^{2})}{2\alpha^{2}-\alpha\hat{\alpha}}e^{-\hat{\alpha}(t-s)+\hat{\varepsilon} |s|}.
\eeaa
Clearly, the proof above is also valid for proving the other terms in the right-hand
side in \x{c44}. Thus we can rewrite the inequality \x{c44} as
\bea\lb{c45}\EM\E \left\|\hat{\Phi}(t,s)\hat{P}(s)\hat{P}(s)-\Phi(t,s)P(s)\hat{P}(s)\right\|^{2}\nn\\
\LE \left(\frac{8M\hat{M}h^{2}}{\alpha-\hat{\alpha}}+
\frac{64M\hat{M}(b^{2}+g^{2}h^{2})}{2\alpha^{2}-\alpha\hat{\alpha}}\right)
e^{-\hat{\alpha}(t-s)+\hat{\varepsilon}|s|}\nn\\
\EQ \frac{320M\tilde{M}}{\alpha-\hat{\alpha}}e^{-\hat{\alpha}(t-s)+\hat{\varepsilon}|s|}.
\eea

On the other hand, since $\hat{P}(s)$ and $\hat{Q}(s)$ are complementary projections for each $s\in I$, it follows from\x{c2}, \x{x3} and \x{c33} that
\bea\lb{c46}
\EM
\E \left\|\hat{\Phi}(t,s)\hat{P}(s)\hat{Q}(s)-\Phi(t,s)P(s)\hat{Q}(s)\right\|^{2}
= \E \left\|\Phi(t,s)P(s)\hat{Q}(s)\right\|^{2}\nn\\
\LE 2\E \left\|\int^{s}_{t_{0}}\Phi(t)\Phi^{-1}(\tau)P(\tau)H(\tau)\hat{\Phi}(\tau, t)\hat{Q}(t)d\omega(\tau)\right\|^{2}\nn\\
 \EM+2\E \left\|\int^{t}_{t_{0}}\Phi(t)\Phi^{-1}(\tau)P(\tau)\tilde{B}(\tau)\hat{\Phi}(\tau, t)\hat{Q}(t)d\tau \right\|^{2}\nn\\
 \LE \frac{40M\tilde{M}}{\alpha-\hat{\alpha}}e^{-\hat{\alpha}(t-s)+\hat{\varepsilon}|s|}.
\eea
Combining \x{c45} and \x{c46} yields
\beaa \EM \E\|\hat{\Phi}(t,s)\hat{P}(s)-\Phi(t,s)P(s)\|^{2}\\
\EQ \E\|\hat{\Phi}(t,s)\hat{P}(s)(\hat{P}(s)+\hat{Q}(s))
-\Phi(t,s)P(s)(\hat{P}(s)+\hat{Q}(s))\|^{2}\\
\LE \frac{720M\tilde{M}}{\alpha-\hat{\alpha}}e^{-\hat{\alpha}(t-s)+\hat{\varepsilon}|s|}.
\eeaa

Similarly, by \x{c18} we obtain
\bea\lb{c48}\EM\E \left\|\hat{\Phi}(t,s)\hat{Q}(s)\hat{Q}(s)-\Phi(t,s)Q(s)\hat{Q}(s)\right\|^{2}\le \frac{320M\tilde{M}}{\alpha-\hat{\alpha}}e^{-\hat{\alpha}(s-t)+\hat{\varepsilon}|s|}.
\eea
On the other hand, since $\hat{P}(s)$ and $\hat{Q}(s)$ are complementary projections for each $s\in I$, by \x{c31} we obtain
\bea\lb{c47}\EM\E \left\|\hat{\Phi}(t,s)\hat{Q}(s)\hat{P}(s)-\Phi(t,s)Q(s)\hat{P}(s)\right\|^{2}\le \frac{40M\tilde{M}}{\alpha-\hat{\alpha}}e^{-\hat{\alpha}(s-t)+\hat{\varepsilon}|s|}.
\eea
Combining \x{c48} and \x{c47} yields
\beaa \EM \E\|\hat{\Phi}(t,s)\hat{Q}(s)-\Phi(t,s)Q(s)\|^{2}\\
\EQ \E\|\hat{\Phi}(t,s)\hat{Q}(s)(\hat{P}(s)+\hat{Q}(s))
-\Phi(t,s)Q(s)(\hat{P}(s)+\hat{Q}(s))\|^{2}\\
\LE \frac{720M\tilde{M}}{\alpha-\hat{\alpha}}e^{-\hat{\alpha}(s-t)+\hat{\varepsilon}|s|}.
\eeaa
This completes the proof of the theorem. \hspace{\stretch{1}}$\Box$
}

\begin{remark}\lb{rem32}
  Since $I=[t_{0},+\oo)$, the second-moment Lyapunov exponent is bounded by  $-\hat{\alpha}$ for any fixed $b, h >0$, i.e.,
  \[\lim_{t\rightarrow +\oo}\frac{1}{t}\log \E\|\hat{\Phi}(t,s)\hat{P}(s)-\Phi(t,s)P(s)\|^{2}=-\hat{\alpha}<0.\]
  This shows that in the stable direction,  any two solutions $\hat{\Phi}(t,s)\hat{P}(s)\xi$ and $\Phi(t,s)P(s)\xi$ with the same initial condition are forward asymptotic in the mean-square sense. Furthermore, since $M=8b^{2}+8g^{2}h^{2}+\alpha h^2$, for each fixed $T_{1}\in (s,+\oo)$ and $T_{2} \in (t_{0},s)$, we have
  \[\lim_{b, h\rightarrow 0}\sup_{t\in [s,T_{1}]} \E\|\hat{\Phi}(t,s)\hat{P}(s)-\Phi(t,s)P(s)\|^{2}=0,\]
  and
    \[\lim_{b, h\rightarrow 0}\sup_{t\in [T_{2},s]} \E\|\hat{\Phi}(t,s)\hat{Q}(s)-\Phi(t,s)Q(s)\|^{2}=0.\]
  This means that the solution $\hat{\Phi}(t,s)\hat{P}(s)$ (or $\hat{\Phi}(t,s)\hat{Q}(s)$) of the perturbed system {\rm\x{a3}} approaches uniformly the solution $\Phi(t,s)P(s)$ (or $\Phi(t,s)Q(s)$) of the system {\rm \x{a1}} in the mean-square sense on any compact interval.
\end{remark}

\section{Robustness of NMS-ED on the half line $(-\oo,t_{0}]$}
\setcounter{equation}{0} \noindent

In this section  we deal with the  robustness of NMS-ED on $I= (-\oo,t_{0}]$, which is analogous to the case
$[t_{0},+\oo)$. So in what follows, we highlight the main steps of the proof which only indicate the major differences.

\begin{theorem}\lb{main31}
  The assertion in Theorem {\rm\ref{main23}} remains true for $I= (-\oo,t_{0}]$.
\end{theorem}

\mbox{\bf Proof of Theorem \ref{main31}.}
Consider the Banach space
\be\lb{d1}\mathscr{L}_{d}:=\{\hat{\Phi}: I^{2}_{\leq} \rightarrow
\mathfrak{B}(L^{2}(\Omega, \R^{n})): ~\hat{\Phi} {\rm~is ~continuous ~and~}
\|\hat{\Phi}\|_{d}<\oo\}\ee with the norm \be\lb{d2}
\|\hat{\Phi}\|_{d}=\sup\left\{(\E\|\hat{\Phi}(t, s)\|^2)^{\frac{1}{2}}
e^{-\frac{\varepsilon}{2} |s|}:
(t,s)\in I^{2}_{\leq}\right\}.\ee

Following the same steps as in the proof of Theorem 1, we establish the following statements.

\begin{lemma}\lb{lem41}
  For  any given initial value $\xi_{0} \in \R^{n}$, {\rm \x{a3}} has a unique solution $V(t,s) \xi_{0}$ with $V \in (\mathscr{L}_{d},
\|\cdot\|_{d})$ such that
   \bea\lb{d10} \EM V(t,s)=
\Phi(t)\Phi^{-1}(s)Q(s)- \int^{s}_{t}\Phi(t)\Phi^{-1}(\tau)Q(\tau)H(\tau)V(\tau, s)d\omega(\tau)
\nn\\
\EM - \int^{s}_{t}\Phi(t)\Phi^{-1}(\tau)Q(\tau)\tilde{B}(\tau)V(\tau, s)d\tau +\int^{t}_{-\oo}\Phi(t)\Phi^{-1}(\tau)P(\tau)H(\tau)V(\tau, s)d\omega(\tau)\nn\\
\EM +\int^{t}_{-\oo}\Phi(t)\Phi^{-1}(\tau)P(\tau)\tilde{B}(\tau)V(\tau, s)d\tau.
\eea
\end{lemma}

\begin{lemma}\lb{lem42}
  For  any $u\in (t,s)$ in $I$, we have
   \[V(s,t)=V(s,u)V(u,t)\]
in the sense of  $(\mathscr{L}_{d},
\|\cdot\|_{d})$.
\end{lemma}

\begin{lemma}\lb{lem43}Given $s\in I$, if $y(t):=\tilde{\Lambda}(t,s) \xi:(-\oo,s]\rightarrow L^{2}(\Omega, \R^{n})$ is a solution of {\rm\x{a3}} with $y(s)= \tilde{\Lambda}(s,s) \xi= \xi$ such that $\tilde{\Lambda}$ is bounded in  $(\mathscr{L}_{d},
\|\cdot\|_{d})$.
Then
   \beaa \EM y(t)=
\Phi(t)\Phi^{-1}(s)Q(s)\xi- \int^{s}_{t}\Phi(t)\Phi^{-1}(\tau)Q(\tau)H(\tau)y(\tau)d\omega(\tau)
\nn\\
\EM - \int^{s}_{t}\Phi(t)\Phi^{-1}(\tau)Q(\tau)\tilde{B}(\tau)y(\tau)d\tau +\int^{t}_{-\oo}\Phi(t)\Phi^{-1}(\tau)P(\tau)H(\tau)y(\tau)d\omega(\tau)\nn\\
\EM +\int^{t}_{-\oo}\Phi(t)\Phi^{-1}(\tau)P(\tau)\tilde{B}(\tau)y(\tau)d\tau.
\eeaa
\end{lemma}
For each $t\in I$, define linear operators as
\be\lb{d11} \hat{Q}(t)=\hat{\Phi}(t,t_{0})V(t_{0},t_{0})\hat{\Phi}(t_{0},t)\quad {\rm and}\quad
\hat{P}(t)=Id-\hat{Q}(t),
\ee
where $t_{0}$ is the right boundary point of the interval $I$.
\begin{lemma}\lb{lem44}
  The operator $\hat{P}(t)$ are linear projections for $t\in I$, and {\rm \x{x4}} holds.
\end{lemma}

\begin{lemma}\lb{lem45}
 For  any given initial value $\xi_{0} \in \R^{n}$, the function  $\hat{Q}(t)\hat{\Phi}(t,s)\xi_{0}$ is a solution of  {\rm \x{a3}} with $\hat{Q}(t)\hat{\Phi}(t,s)$ is bounded  in  $(\mathscr{L}_{d},
\|\cdot\|_{d})$.
\end{lemma}

\begin{lemma}\lb{lem46}
 For  any given initial value $\xi_{0} \in \R^{n}$, the function  $\hat{Q}(t)\hat{\Phi}(t,s)\xi_{0}$ is a solution of  {\rm \x{a3}} with  $(t,s)\in I^{2}_{\leq}$ such that
\bea\lb{d3} \EM \hat{\Phi}(t,s)\hat{Q}(s)=
\Phi(t)\Phi^{-1}(s)Q(s)\hat{Q}(s)- \int^{s}_{t}\Phi(t)\Phi^{-1}(\tau)Q(\tau)H(\tau)\hat{\Phi}(\tau, s)\hat{Q}(s)d\omega(\tau)
\nn\\
\EM - \int^{s}_{t}\Phi(t)\Phi^{-1}(\tau)Q(\tau)\tilde{B}(\tau)\hat{\Phi}(\tau, s)\hat{Q}(s)d\tau +\int^{t}_{-\oo}\Phi(t)\Phi^{-1}(\tau)P(\tau)H(\tau)\hat{\Phi}(\tau, s)\hat{Q}(s)d\omega(\tau)\nn\\
\EM +\int^{t}_{-\oo}\Phi(t)\Phi^{-1}(\tau)P(\tau)\tilde{B}(\tau)\hat{\Phi}(\tau, s)\hat{Q}(s)d\tau.
\eea
\end{lemma}

\begin{lemma}\lb{lem47}
 For  any given initial value $\xi_{0} \in \R^{n}$, the function  $\hat{P}(t)\hat{\Phi}(t,s)\xi_{0}$ is a solution of  {\rm \x{a3}} with  $(t,s)\in I^{2}_{\geq}$ such that
\bea\lb{d4} \EM \hat{\Phi}(t,s)\hat{P}(s)=
\Phi(t)\Phi^{-1}(s)P(s)\hat{P}(s)- \int^{t_{0}}_{t}\Phi(t)\Phi^{-1}(\tau)Q(\tau)H(\tau)\hat{\Phi}(\tau, s)\hat{Q}(s)d\omega(\tau)
\nn\\
\EM - \int^{t_{0}}_{t}\Phi(t)\Phi^{-1}(\tau)Q(\tau)\tilde{B}(\tau)\hat{\Phi}(\tau, s)\hat{Q}(s)d\tau  +\int^{t}_{s}\Phi(t)\Phi^{-1}(\tau)P(\tau)H(\tau)\hat{\Phi}(\tau, s)\hat{P}(s)d\omega(\tau)\nn\\
\EM +\int^{t}_{s}\Phi(t)\Phi^{-1}(\tau)P(\tau)\tilde{B}(\tau)\hat{\Phi}(\tau, s)\hat{P}(s)d\tau.
\eea
\end{lemma}

Proceed as in the proof of Theorem \ref{main23}.
 Squaring both sides of  \x{d3}, and taking expectations, we obtain
\be\lb{d5}\E\|\hat{\Phi}(t,s)\hat{Q}(s)\|^{2}\leq 5M
e^{(-\frac{\alpha}{2}+\frac{10M\tilde{M}}{\alpha})(s-t)+\varepsilon |s|}\E\|\hat{Q}(s)\|^{2},\quad \forall~ (t,s)\in
I^{2}_{\leq}.\ee
Similarly, Squaring both sides of  \x{d4}, and taking expectations, we obtain
\be\lb{d6}\E\|\hat{\Phi}(t,s)\hat{P}(s)\|^{2}\leq 5M
e^{(-\frac{\alpha}{2}+\frac{10M\tilde{M}}{\alpha})(t-s)+\varepsilon |s|}\E\|\hat{P}(s)\|^{2},\quad \forall~ (t,s)\in
I^{2}_{\geq}.\ee
Meanwhile, multiplying \x{d3} with $P(t)$ and  \x{d4} with $Q(t)$ on the left side, respectively, and let $t=s$, we obtain
\[ \E\|P(t)\hat{Q}(t)\|^{2}\le \frac{10M^{2}\tilde{M}}{\alpha(\alpha+\tilde{\alpha}-\varepsilon)}\E\|\hat{Q}(t)\|^{2},
\]
and
\[ \E\|Q(t)\hat{P}(t)\|^{2}\le \frac{10M^{2}\tilde{M}}{\alpha(\alpha+\tilde{\alpha}-\varepsilon)}\E\|\hat{P}(t)\|^{2}.
\]
Since \[ \E\|P(t)\|^{2}\le M e^{\varepsilon |t|},\quad  \E\|Q(t)\|^{2}\le M e^{\varepsilon |t|},\]
and
$\hat{P}(t)-P(t)=Q(t)\hat{P}(t)-P(t)\hat{Q}(t)$, for sufficiently small $b$ and $h$, we obtain the bounds for the projections $\hat{P}(t)$ and $\hat{Q}(t)$ as follows:
\be\lb{d9} \E\|\hat{P}(t)\|^{2}\le 8M e^{\varepsilon |t|}\quad {\rm and}\quad \E\|\hat{Q}(t)\|^{2}\le 8M e^{\varepsilon |t|}.\ee
By \x{d5}, \x{d6}, using \x{d9} we obtain
\[\E\|\hat{\Phi}(t,s)\hat{P}(s)\|^{2}\leq 40M^{2}
e^{(-\frac{\alpha}{2}+\frac{10M\tilde{M}}{\alpha})(t-s)+2\varepsilon |s|},\quad \forall~ (t,s)\in
I^{2}_{\geq},\]
and
\[\E\|\hat{\Phi}(t,s)\hat{Q}(s)\|^{2}\leq 40M
e^{(-\frac{\alpha}{2}+\frac{10M\tilde{M}}{\alpha})(s-t)+2\varepsilon |s|},\quad \forall~ (t,s)\in
I^{2}_{\leq}.\]
This completes the proof of the theorem. \hspace{\stretch{1}}$\Box$

\section{Robustness of NMS-ED on the whole $\R$}
\setcounter{equation}{0} \noindent
In this section  we consider   the robustness of  NMS-ED  on the whole $I= \R$.  From the last two sections we know that if \x{x1} holds, the perturbed equation \x{a3} remains NMS-ED on $[t_{0},+\oo)$ with the operators:
\beaa \hat{P}_{+}(t)=\hat{\Phi}(t,t_{0})U(t_{0},t_{0})\hat{\Phi}(t_{0},t),\quad
\hat{Q}_{+}(t)=Id-\hat{P}_{+}(t),
\eeaa
and on $(-\oo, t_{0}]$ with the operators:
\beaa\hat{Q}_{-}(t)=\hat{\Phi}(t,t_{0})V(t_{0},t_{0})\hat{\Phi}(t_{0},t),\quad
\hat{P}_{-}(t)=Id-\hat{Q}_{-}(t).
\eeaa
The most important part in this section is to show that  \x{a3} has an NMS-ED on both half lines with the same projections. For this purpose we introduce modified  projections, which combines the advantages of projections $\hat{P}_{+}(t)$ and $\hat{Q}_{-}(t)$. Actually, this technique has been used in a lot of papers to deal with this problem, see e.g., \cite{bsv09, bv08, pa-84, pa-88, pop-06, pop-09} for details.

In the following, for convenience and brevity, let us denote by $G(t,s)$ the Green function of \x{a1}:
\[G(t,s):=\left\{ \begin{aligned}& P(t)\Phi(t,s), \quad & \forall~ (t,s)\in
\R^{2}_{\geq}, \\
&       -Q(t)\Phi(t,s), \quad  & \forall~ (t,s)\in
\R^{2}_{\leq}.\end{aligned}
\right.\]
Green function is a classical concept in the study of exponential dichotomy as for example \cite{chi, en}. Now we deal with the  robustness of NMS-ED for \x{a1} on the whole $\R$.
\begin{theorem}\lb{main51}
  The assertion in Theorem {\rm\ref{main23}} remains true for $I= \R$.
\end{theorem}
\mbox{\bf Proof of Theorem \ref{main51}.}
Consider the Banach spaces
\beaa\mathscr{L}_{c}\EQ \{\hat{\Phi}: \R^{2}_{\geq} \rightarrow
\mathfrak{B}(L^{2}(\Omega, \R^{n})): ~\hat{\Phi} {\rm~is ~continuous ~and~}
\|\hat{\Phi}\|_{c}<\oo\},\eeaa
and
\beaa \mathscr{L}_{d}\EQ \{\hat{\Phi}: \R^{2}_{\leq} \rightarrow
\mathfrak{B}(L^{2}(\Omega, \R^{n})): ~\hat{\Phi} {\rm~is ~continuous ~and~}
\|\hat{\Phi}\|_{d}<\oo\}
\eeaa with the norm
\beaa\|\hat{\Phi}\|_{c}\EQ\sup_{(t,s)\in \R^{2}_{\geq}}\left\{(\E\|\hat{\Phi}(t, s)\|^2)^{\frac{1}{2}}
e^{-\frac{\varepsilon}{2} |s|}
\right\},\eeaa
and
\beaa\|\hat{\Phi}\|_{d}\EQ\sup_{(t,s)\in \R^{2}_{\leq}}\left\{(\E\|\hat{\Phi}(t, s)\|^2)^{\frac{1}{2}}
e^{-\frac{\varepsilon}{2} |s|}
\right\}\eeaa
respectively. Define operator $\Gamma_{1}:\mathscr{L}_{c} \rightarrow \mathscr{L}_{c}$ by
\beaa (\Gamma_{1} U)(t,s)\EQ \Phi(t)\Phi^{-1}(s)P(s)+ \int^{\oo}_{s}G(t,\tau)H(\tau)U(\tau, s)d\omega(\tau)
\nn\\
\EM + \int^{\oo}_{s}G(t,\tau)\tilde{B}(\tau)U(\tau, s)d\tau, \eeaa
and operator $\Gamma_{2}:\mathscr{L}_{d} \rightarrow \mathscr{L}_{d}$,
\beaa (\Gamma_{2} V)(t,s)\EQ \Phi(t)\Phi^{-1}(s)Q(s)+ \int^{s}_{-\oo}G(t,\tau)H(\tau)V(\tau, s)d\omega(\tau)
\nn\\
\EM + \int^{s}_{-\oo}G(t,\tau)\tilde{B}(\tau)V(\tau, s)d\tau. \eeaa
Similar arguments to those in the proofs of Lemma \ref{lem32} and Lemma \ref{lem41} can be used to deduce that
\beaa\|\Gamma_{1} U_{1}-\Gamma_{1}  U_{2}\|_{c}\LE \theta \|U_{1}-U_{2}\|_{c},\\
\|\Gamma_{2} V_{1}-\Gamma_{2}  V_{2}\|_{d}\LE \theta \|V_{1}-V_{2}\|_{d}.\eeaa
with $\theta=\sqrt{\frac{10M\tilde{M}}{\alpha^{2}} }<1 $. Thus we have the following lemma.

\begin{lemma}
  Operators $\Gamma_{1}$, $\Gamma_{2}$ have unique fixed points $U\in (\mathscr{L}_{c},
\|\cdot\|_{c})$, respectively $V\in (\mathscr{L}_{d},
\|\cdot\|_{d})$ such that
\end{lemma}
\beaa  U(t,s)\EQ \Phi(t)\Phi^{-1}(s)P(s)+ \int^{\oo}_{s}G(t,\tau)H(\tau)U(\tau, s)d\omega(\tau)
\nn\\
\EM + \int^{\oo}_{s}G(t,\tau)\tilde{B}(\tau)U(\tau, s)d\tau, \eeaa
and
\beaa   V(t,s)\EQ \Phi(t)\Phi^{-1}(s)Q(s)+ \int^{s}_{-\oo}G(t,\tau)H(\tau)V(\tau, s)d\omega(\tau)
\nn\\
\EM + \int^{s}_{-\oo}G(t,\tau)\tilde{B}(\tau)V(\tau, s)d\tau. \eeaa

Repeating arguments in the proofs of Theorem \ref{main23} and Theorem \ref{main31}  we obtain the following statements.

\begin{lemma}
  For  any $u\in (s,t)$ in $I$, we have
   \[U(t,s)=U(t,u)U(u,s)\]
in the sense of  $(\mathscr{L}_{c},
\|\cdot\|_{c})$, respectively,  \[V(t,s)=V(t,u)V(u,s)\]
in the sense of  $(\mathscr{L}_{d},
\|\cdot\|_{d})$.
\end{lemma}

\begin{lemma}Given $s\in I$, if $x(t)=\Lambda(t,s) \xi:[s,+\oo)\rightarrow L^{2}(\Omega, \R^{n})$ (respectively,  $y(t):=\tilde{\Lambda}(t,s) \xi:(-\oo,s]\rightarrow L^{2}(\Omega, \R^{n})$ ) is a solution of {\rm\x{a3}} with $x(s)= \Lambda(s,s) \xi= \xi$ (respectively, $y(s)= \tilde{\Lambda}(s,s) \xi= \xi$) such that $\Lambda$ (respectively, $\tilde{\Lambda}$) is bounded in  $(\mathscr{L}_{c},
\|\cdot\|_{c})$ (respectively, $(\mathscr{L}_{d},
\|\cdot\|_{d})$).
Then
\bea\lb{f1}  x(t)\EQ
\Phi(t)\Phi^{-1}(s)P(s)\xi+ \int^{\oo}_{s}G(t,\tau)H(\tau)x(\tau)d\omega(\tau)
\nn\\
\EM + \int^{\oo}_{s}G(t,\tau)\tilde{B}(\tau)x(\tau)d\tau,
\eea
and
   \bea\lb{f2} y(t)\EQ
\Phi(t)\Phi^{-1}(s)Q(s)\xi+ \int^{s}_{-\oo}G(t,\tau)H(\tau)y(\tau)d\omega(\tau)
\nn\\
\EM +\int^{s}_{-\oo}G(t,\tau)\tilde{B}(\tau)y(\tau)d\tau.
\eea
\end{lemma}

Now we present that projection $S=\hat{P}_{+}(t_{0})+\hat{Q}_{-}(t_{0})$ is invertible for some $t_{0}\in \R$ with $b$ and $h$ are sufficiently small. Using this result, we are able to define modified operators.

\begin{lemma}
  If $b$ and $h$ are sufficiently small, then the operator $S=\hat{P}_{+}(t_{0})+\hat{Q}_{-}(t_{0})$ is invertible.
\end{lemma}
\prf{We first derive $\hat{P}_{+}(t_{0})P(t_{0})=\hat{P}_{+}(t_{0})$. In fact, following the same procedure as we did for Lemma \ref{lem33} we find that
\be\lb{f3}U(t,s)=U(t,s)P(s).\ee
Since $\hat{P}_{+}(t_{0})=U(t_{0},t_{0})$, by \x{f3} with $t=s=t_{0}$ we have
\be\lb{f4}\hat{P}_{+}(t_{0})P(t_{0})=\hat{P}_{+}(t_{0}).\ee
In addition, we have (see \x{c22})
\be\lb{f5}P(t_{0})\hat{P}_{+}(t_{0})=P(t_{0}).\ee
Since $\hat{Q}_{-}(t_{0})=V(t_{0},t_{0})$, a similar argument using Lemma \ref{lem42} with $t=s=t_{0}$ yields
\be\lb{f6}\hat{Q}_{-}(t_{0})Q(t_{0})=\hat{Q}_{-}(t_{0}).\ee
On the other hand, it follows from $\hat{Q}_{-}(t)=\hat{\Phi}(t,t_{0})V(t_{0},t_{0})\hat{\Phi}(t_{0},t)$ and \x{d10} with $t=s=t_{0}$ that
 \bea\lb{f7} \hat{Q}_{-}(t_{0}) \EQ V(t_{0},t_{0})=
Q(t_{0})
 +\int^{t_{0}}_{-\oo}\Phi(t_{0})\Phi^{-1}(\tau)P(\tau)H(\tau)V(\tau, t_{0})d\omega(\tau)\nn\\
\EM +\int^{t_{0}}_{-\oo}\Phi(t_{0})\Phi^{-1}(\tau)P(\tau)\tilde{B}(\tau)V(\tau, t_{0})d\tau.
\eea
Since $P(t_{0})$ and  $Q(t_{0})$ are complementary projections, multiplies \x{f7} on the left with $Q(t_{0})$. This gives
\bea\lb{f8} Q(t_{0})\hat{Q}_{-}(t_{0})=Q(t_{0}).
\eea
We now consider the linear operators
\be\lb{f9} S_{1}:=Id -P(t_{0})+\hat{P}_{+}(t_{0}) \quad {\rm and} \quad
T_{1}:= Id+ P(t_{0})-\hat{P}_{+}(t_{0}).
\ee
It follows easily from \x{f4} and \x{f5} that $S_{1}T_{1}=T_{1}S_{1}=Id$. Therefore, $S_{1}$ is invertible and $S_{1}^{-1}=T_{1}$. In addition, using again \x{f5} we obtain
\bea\lb{f10} S_{1}-Id \EQ \hat{P}_{+}(t_{0})-P(t_{0})\nn\\
\EQ \hat{P}_{+}(t_{0})-P(t_{0})\hat{P}_{+}(t_{0})\nn\\
\EQ Q(t_{0})\hat{P}_{+}(t_{0}).
\eea
By \x{c21}, we have
 \bea\lb{f11} Q(t_{0})\hat{P}_{+}(t_{0}) \EQ
 -\int^{\oo}_{t_{0}}\Phi(t_{0})\Phi^{-1}(\tau)Q(\tau)H(\tau)U(\tau, t_{0})d\omega(\tau)\nn\\
\EM -\int^{\oo}_{t_{0}}\Phi(t_{0})\Phi^{-1}(\tau)Q(\tau)\tilde{B}(\tau)U(\tau, t_{0})d\tau.
\eea
To estimate the bounds  of the integral in the mean square sense, we need to find out the bounds for $U(t, t_{0})$ with $t\ge t_{0}$  Squaring both sides of \x{c4},
taking expectations, and proceeding as in the proof of Theorem \ref{main23}, for any $t\ge t_{0}$,  we have
\bea\lb{f12} \EM \E\|U(t,t_{0})\|^{2} \le 5\E \|\Phi(t)\Phi^{-1}(t_{0})P(t_{0})\|^{2} +
5\E \left\|\int^{t}_{t_{0}}\Phi(t)\Phi^{-1}(\tau)P(\tau)H(\tau)U(\tau, t_{0})d\omega(\tau)\right\|^{2}\nn\\
\EM + 5\E \left\|\int^{t}_{t_{0}}\Phi(t)\Phi^{-1}(\tau)P(\tau)\tilde{B}(\tau)U(\tau, t_{0})d\tau \right\|^{2}+ 5\E \left\|\int^{\oo}_{t}\Phi(t)\Phi^{-1}(\tau)Q(\tau)H(\tau)U(\tau, t_{0})d\omega(\tau) \right\|^{2}\nn\\
\EM + 5\E \left\|\int^{\oo}_{t}\Phi(t)\Phi^{-1}(\tau)Q(\tau)\tilde{B}(\tau)U(\tau, t_{0})d\tau \right\|^{2}\nn\\
\LE 5Me^{-\frac{\alpha}{2}(t-t_{0})+\varepsilon |t_{0}|}+ \frac{5M\tilde{M}}{\alpha}\left(\int^{t}_{t_{0}}
e^{-\frac{\alpha}{2}(t-\tau)}\E\|U(\tau,t_{0})\|^{2} d \tau+
\int^{\oo}_{t}e^{-\frac{\alpha}{2}(\tau-t)}\E\|U(\tau,t_{0})\|^{2} d \tau\right)\nn\\
\LE 5M
e^{-\hat{\alpha}(t-t_{0})+\varepsilon |t_{0}|}.
\eea
By \x{f10}, using \x{f11} and \x{f12}, we obtain
\bea\lb{f13} \E\|S_{1}-Id\|^{2}\EQ \E\|Q(t_{0})\hat{P}_{+}(t_{0})\|^{2}\nn\\
\LE 2\E \left\|\int^{\oo}_{t_{0}}\Phi(t_{0})\Phi^{-1}(\tau)Q(\tau)H(\tau)U(\tau, t_{0})d\omega(\tau)\right\|^{2}\nn\\
\EM+2 \E \left\|\int^{\oo}_{t_{0}}\Phi(t_{0})\Phi^{-1}(\tau)Q(\tau)\tilde{B}(\tau)U(\tau, t_{0})d\tau\right\|^{2}\nn\\
\LE 2\int^{\oo}_{t} \E\|\Phi(t)\Phi^{-1}(\tau)Q(\tau)\|^{2}
\E\|H(\tau)\|^{2} \E\|U(\tau, t_{0})\|^{2} d\tau \nn\\
\EM + 2 \left(\int^{\oo}_{t_{0}} \E \|\Phi(t_{0})\Phi^{-1}(\tau)Q(\tau)\|\E\|\tilde{B}(\tau)\|^{\frac{1}{2}}d \tau\right)\nn\\
 \EM\times \left(\int^{\oo}_{t_{0}} \E \|\Phi(t_{0})\Phi^{-1}(\tau)Q(\tau)\|\E\|\tilde{B}(\tau)\|^{\frac{3}{2}} \E\|U(\tau, t_{0})\|^{2}d \tau\right)\nn\\
 \LE \frac{10M^{2}\tilde{M}}{\alpha}\int^{\oo}_{t_{0}}
 e^{-(\alpha+\tilde{\alpha}-\varepsilon)(\tau-t_{0})}d \tau\nn\\
 \LE \frac{10M^{2}\tilde{M}}{\alpha(\alpha+\tilde{\alpha}-\varepsilon)}.
\eea

Meanwhile, we consider the linear operators
\be\lb{f14} S_{2}:=Id -Q(t_{0})+\hat{Q}_{-}(t_{0}) \quad {\rm and} \quad
T_{2}:= Id+ Q(t_{0})-\hat{Q}_{-}(t_{0}).
\ee
It follows easily from \x{f6} and \x{f8} that $S_{2}T_{2}=T_{2}S_{2}=Id$. Therefore, $S_{2}$ is invertible and $S_{2}^{-1}=T_{2}$. In addition, using again \x{f8} we obtain
\bea\lb{f15} S_{2}-Id \EQ \hat{Q}_{-}(t_{0})-Q(t_{0})\nn\\
\EQ \hat{Q}_{-}(t_{0})-Q(t_{0})\hat{Q}_{-}(t_{0})\nn\\
\EQ P(t_{0})\hat{Q}_{-}(t_{0}).
\eea
By \x{f7},
\bea\lb{f16} P(t_{0})\hat{Q}_{-}(t_{0})\EQ
\int^{t_{0}}_{-\oo}\Phi(t_{0})\Phi^{-1}(\tau)P(\tau)H(\tau)V(\tau, t_{0})d\omega(\tau)\nn\\
\EM +\int^{t_{0}}_{-\oo}\Phi(t_{0})\Phi^{-1}(\tau)P(\tau)\tilde{B}(\tau)V(\tau, t_{0})d\tau.
\eea
Similarly, for any $t\le t_{0}$, one can deduce from \x{d10} that
\bea\lb{f17} \EM \E\|V(t,t_{0})\|^{2} \le 5\E \|\Phi(t)\Phi^{-1}(t_{0})Q(t_{0})\|^{2} +
5\E \left\|\int^{t_{0}}_{t}\Phi(t)\Phi^{-1}(\tau)Q(\tau)H(\tau)V(\tau, s)d\omega(\tau)\right\|^{2}\nn\\
\EM + 5\E \left\|\int^{t_{0}}_{t}\Phi(t)\Phi^{-1}(\tau)Q(\tau)\tilde{B}(\tau)U(\tau, s)d\tau \right\|^{2}+ 5\E \left\|\int^{t}_{-\oo}\Phi(t)\Phi^{-1}(\tau)P(\tau)H(\tau)V(\tau, t_{0})d\omega(\tau) \right\|^{2}\nn\\
\EM + 5\E \left\|\int^{t}_{-\oo}\Phi(t)\Phi^{-1}(\tau)P(\tau)\tilde{B}(\tau)V(\tau, t_{0})d\tau \right\|^{2}\nn\\
\LE 5Me^{-\frac{\alpha}{2}(t_{0}-t)+\varepsilon |t_{0}|}+ \frac{5M\tilde{M}}{\alpha}\left(\int^{t_{0}}_{t}
e^{-\frac{\alpha}{2}(\tau-t)}\E\|V(\tau,t_{0})\|^{2} d \tau+
\int^{t}_{-\oo}e^{-\frac{\alpha}{2}(t-\tau)}\E\|V(\tau,t_{0})\|^{2} d \tau\right)\nn\\
\LE 5M
e^{-\hat{\alpha}(t_{0}-t)+\varepsilon |t_{0}|}.
\eea
Therefore, by \x{f15}, using \x{f16} and \x{f17} we obtain
\bea\lb{f18} \E\|S_{2}-Id\|^{2}\EQ \E\|P(t_{0})\hat{Q}_{-}(t_{0})\|^{2}
\le \frac{10M^{2}\tilde{M}}{\alpha(\alpha+\tilde{\alpha}-\varepsilon)}.
\eea

On the other side, it follows easily from \x{f8} that $P(t_{0})\hat{P}_{-}(t_{0})=\hat{P}_{-}(t_{0})$. Using also \x{f5} yields
\beaa \hat{P}_{+}(t_{0})+\hat{Q}_{-}(t_{0})-Id \EQ \hat{P}_{+}(t_{0})
-P(t_{0})+P(t_{0})-\hat{P}_{-}(t_{0})\\
\EQ \hat{P}_{+}(t_{0})-P(t_{0})\hat{P}_{+}(t_{0})+P(t_{0})-P(t_{0})\hat{P}_{-}(t_{0})\\
\EQ Q(t_{0})\hat{P}_{+}(t_{0}) + P(t_{0})\hat{Q}_{-}(t_{0}).
\eeaa
By \x{f13} and \x{f18} we obtain
\bea\lb{f19} \E\|\hat{P}_{+}(t_{0})+\hat{Q}_{-}(t_{0})-Id\|^{2}
\EQ \E\|Q(t_{0})\hat{P}_{+}(t_{0}) + P(t_{0})\hat{Q}_{-}(t_{0})\|^{2}\nn\\
\LE 2\E\|Q(t_{0})\hat{P}_{+}(t_{0})\|^{2}+2\E\|P(t_{0})\hat{Q}_{-}(t_{0})\|^{2}\nn\\
\LE \frac{20M^{2}\tilde{M}}{\alpha(\alpha+\tilde{\alpha}-\varepsilon)}.
\eea
Moreover,
\bea\lb{f20} S \EQ \hat{P}_{+}(t_{0})+\hat{Q}_{-}(t_{0})\nn\\
\EQ (\hat{P}_{+}(t_{0})+Q(t_{0})) +(P(t_{0})+\hat{Q}_{-}(t_{0}))-Id\nn\\
\EQ S_{1}+S_{2}-Id.
\eea
Since $\tilde{M}:=8b^{2}+8g^{2}h^{2}+\alpha h^2$, by \x{f13}, respectively, \x{f18}, we can make invertible operator $S_{1}$ and $S_{2}$ such that $\E\|S_{1}-Id\|^{2}$ and  $\E\|S_{2}-Id\|^{2}$ as small as desired with $b$ and $h$ sufficiently small.
So if taking $b$ and $h$ sufficiently small, it follows from \x{f19} and \x{f20} that $S=\hat{P}_{+}(t_{0})+\hat{Q}_{-}(t_{0})$ is invertible.  \hspace{\stretch{1}}$\Box$}

For each $t\in I$, define linear operators as
\be\lb{f21} \tilde{P}(t)=\hat{\Phi}(t,t_{0})SP(t_{0})S^{-1}\hat{\Phi}(t_{0},t)\quad {\rm and}\quad
\tilde{Q}(t)=Id-\tilde{P}(t).
\ee

\begin{lemma}
  The operator $\tilde{P}(t)$ are linear projections for $t\in I$, and {\rm \x{x4}} holds for any $t, s\in \R$.
\end{lemma}
\prf{Obviously,
\beaa \tilde{P}(t)\tilde{P}(t)=\hat{\Phi}(t,t_{0})SP^{2}(t_{0})S^{-1}\hat{\Phi}(t_{0},t)
=\tilde{P}(t).
\eeaa
 Moreover, for any $t, s\in \R$, we obtain
\beaa \tilde{P}(t)\hat{\Phi}(t,s)\EQ \hat{\Phi}(t,t_{0})SP(t_{0})S^{-1}\hat{\Phi}(t_{0},t)\hat{\Phi}(t,s)\\
\EQ \hat{\Phi}(t,s)\hat{\Phi}(s,t_{0})SP(t_{0})S^{-1}\hat{\Phi}(t_{0},s)\\
\EQ \hat{\Phi}(t,s)\tilde{P}(s),
\eeaa
and this completes the proof of the lemma. \hspace{\stretch{1}}$\Box$
}

\begin{lemma}\lb{lem56}
 For  any given initial value $\xi_{0} \in \R^{n}$, the function  $\tilde{P}(t)\hat{\Phi}(t,s)\xi_{0}$ is a solution of  {\rm \x{a3}} with $\tilde{P}(t)\hat{\Phi}(t,s)$ is bounded  in  $(\mathscr{L}_{c},
\|\cdot\|_{c})$, respectively, the function  $\tilde{Q}(t)\hat{\Phi}(t,s)\xi_{0}$ is a solution of  {\rm \x{a3}} with $\tilde{Q}(t)\hat{\Phi}(t,s)$ is bounded  in  $(\mathscr{L}_{d},
\|\cdot\|_{d})$.
\end{lemma}
\prf{In view of \x{f4} and \x{f6}, we have
\beaa SP(t_{0})\EQ\hat{P}_{+}(t_{0})P(t_{0})+\hat{Q}_{-}(t_{0})P(t_{0})=\hat{P}_{+}(t_{0}),\\
SQ(t_{0})\EQ\hat{P}_{+}(t_{0})Q(t_{0})+\hat{Q}_{-}(t_{0})Q(t_{0})=\hat{Q}_{-}(t_{0}).
\eeaa
Thus,
\beaa\tilde{P}(t)\hat{\Phi}(t,s)\EQ\hat{\Phi}(t,t_{0})SP(t_{0})S^{-1}
\hat{\Phi}(t_{0},t)\hat{\Phi}(t,s)= \hat{\Phi}(t,t_{0})\hat{P}_{+}(t_{0})S^{-1}
\hat{\Phi}(t_{0},s),\\
\tilde{Q}(t)\hat{\Phi}(t,s)\EQ\hat{\Phi}(t,t_{0})SQ(t_{0})S^{-1}
\hat{\Phi}(t_{0},t)\hat{\Phi}(t,s)= \hat{\Phi}(t,t_{0})\hat{Q}_{-}(t_{0})S^{-1}
\hat{\Phi}(t_{0},s).
\eeaa
Therefore, it follows from Lemma \ref{lem36} that $\tilde{P}(t)\hat{\Phi}(t,s)\xi_{0}=\hat{\Phi}(t,t_{0})\hat{P}_{+}(t_{0})S^{-1}
\hat{\Phi}(t_{0},s)\xi_{0}$ is a solution of \x{a3} with initial value $S^{-1}
\hat{\Phi}(t_{0},s)\xi_{0} \in\R^{n}$ with $\tilde{P}(t)\hat{\Phi}(t,s)$ is bounded  in  $(\mathscr{L}_{c},
\|\cdot\|_{c})$. Similarly, by Lemma \ref{lem45}, we have $\tilde{Q}(t)\hat{\Phi}(t,s)\xi_{0}$ is a solution of \x{a3} with initial value $S^{-1}
\hat{\Phi}(t_{0},s)\xi_{0} \in\R^{n}$ with $\tilde{Q}(t)\hat{\Phi}(t,s)$ is bounded  in  $(\mathscr{L}_{d},
\|\cdot\|_{d})$. \hspace{\stretch{1}}$\Box$
}

\begin{lemma}
 For  any given initial value $\xi_{0} \in \R^{n}$, the function  $\tilde{P}(t)\hat{\Phi}(t,s)\xi_{0}$ is a solution of  {\rm \x{a3}} with  $(t,s)\in \R^{2}_{\geq}$ such that
\bea\lb{f22} \hat{\Phi}(t,s)\tilde{P}(s) \EQ \Phi(t)\Phi^{-1}(s)P(s)\tilde{P}(s)
+ \int^{\oo}_{s}G(t,\tau)H(\tau)\hat{\Phi}(\tau,s)\tilde{P}(s)d\omega(\tau)\nn\\
\EM+\int^{\oo}_{s}G(t,\tau)\tilde{B}(\tau)\hat{\Phi}(\tau,s)\tilde{P}(s)d\tau,
\eea
and the function  $\tilde{Q}(t)\hat{\Phi}(t,s)\xi_{0}$ is a solution of  {\rm \x{a3}} with  $(t,s)\in \R^{2}_{\leq}$ such that
\bea\lb{f23} \hat{\Phi}(t,s)\tilde{Q}(s) \EQ \Phi(t)\Phi^{-1}(s)Q(s)\tilde{Q}(s)
+\int^{s}_{-\oo}G(t,\tau)H(\tau)\hat{\Phi}(\tau,s)\tilde{Q}(s)d\omega(\tau)\nn\\
\EM+\int^{s}_{-\oo}G(t,\tau)\tilde{B}(\tau)\hat{\Phi}(\tau,s)\tilde{Q}(s)d\tau.
\eea
\end{lemma}
\prf{Let $x(t)=\tilde{P}(t)\hat{\Phi}(t,s)\xi_{0}$ (respectively, $y(t)=\tilde{Q}(t)\hat{\Phi}(t,s)\xi_{0}$) with given $s\in \R$, and denote $\xi=\tilde{P}(s)\xi_{0}$ the initial condition at time $s$.  Clearly, $x(t)$ (respectively, $y(t)$) is a solution of \x{a3} with $x(s)=\tilde{P}(s)\xi=\tilde{P}(s)\tilde{P}(s)\xi_{0}=\xi$ (respectively, $y(s)=\tilde{Q}(s)\xi=\tilde{Q}(s)\tilde{Q}(s)\xi_{0}=\xi$). By Lemma \ref{lem56}, $\tilde{P}(t)\hat{\Phi}(t,s)$ (respectively, $\tilde{Q}(t)\hat{\Phi}(t,s)$) is bounded  in  $(\mathscr{L}_{c},
\|\cdot\|_{c})$ (respectively, $(\mathscr{L}_{d},
\|\cdot\|_{d}))$. Since $\xi_{0}$ is arbitrary in $\R^{n}$, the identity \x{f22} (respectively, \x{f23}) follows now readily from \x{f1} (respectively, \x{f2}). \hspace{\stretch{1}}$\Box$}

Proceed as in the proof of Theorem \ref{main23}.
 Squaring both sides of  \x{f22}, and taking expectations, we obtain
\be\lb{f24}\E\|\hat{\Phi}(t,s)\tilde{P}(s)\|^{2}\leq 5M
e^{(-\frac{\alpha}{2}+\frac{10M\tilde{M}}{\alpha})(s-t)+\varepsilon |s|}\E\|\tilde{P}(s)\|^{2},\quad \forall~ (t,s)\in
\R^{2}_{\geq}.\ee
Similarly, Squaring both sides of  \x{f23}, and taking expectations, we obtain
\be\lb{f25}\E\|\hat{\Phi}(t,s)\tilde{Q}(s)\|^{2}\leq 5M
e^{(-\frac{\alpha}{2}+\frac{10M\tilde{M}}{\alpha})(t-s)+\varepsilon |s|}\E\|\tilde{Q}(s)\|^{2},\quad \forall~ (t,s)\in
\R^{2}_{\leq}.\ee
Meanwhile, multiplying \x{f22} with $Q(t)$ and  \x{f23} with $P(t)$ on the left side, respectively, and let $t=s$, we obtain
\[\E\|Q(t)\tilde{P}(t)\|^{2}\le \frac{10M^{2}\tilde{M}}{\alpha(\alpha+\tilde{\alpha}-\varepsilon)}\E\|\tilde{P}(t)\|^{2},
\]
and
\[ \E\|P(t)\tilde{Q}(t)\|^{2}\le \frac{10M^{2}\tilde{M}}{\alpha(\alpha+\tilde{\alpha}-\varepsilon)}\E\|\tilde{Q}(t)\|^{2}.
\]
Since \[ \E\|P(t)\|^{2}\le M e^{\varepsilon |t|},\quad  \E\|Q(t)\|^{2}\le M e^{\varepsilon |t|},\]
and
$\tilde{P}(t)-P(t)=Q(t)\tilde{P}(t)-P(t)\tilde{Q}(t)$, for sufficiently small $b$ and $h$, we obtain the bounds for the projections $\tilde{P}(t)$ and $\tilde{Q}(t)$ as follows:
\be\lb{f26} \E\|\tilde{P}(t)\|^{2}\le 8M e^{\varepsilon |t|}\quad {\rm and}\quad \E\|\tilde{Q}(t)\|^{2}\le 8M e^{\varepsilon |t|}.\ee
By \x{f24}, \x{f25}, using \x{f26} we obtain
\[\E\|\hat{\Phi}(t,s)\tilde{P}(s)\|^{2}\leq 40M^{2}
e^{(-\frac{\alpha}{2}+\frac{10M\tilde{M}}{\alpha})(t-s)+2\varepsilon |s|},\quad \forall~ (t,s)\in
\R^{2}_{\geq},\]
and
\[\E\|\hat{\Phi}(t,s)\tilde{Q}(s)\|^{2}\leq 40M
e^{(-\frac{\alpha}{2}+\frac{10M\tilde{M}}{\alpha})(s-t)+2\varepsilon |s|},\quad \forall~ (t,s)\in
\R^{2}_{\leq}.\]
This completes the proof of the theorem. \hspace{\stretch{1}}$\Box$

\begin{remark}\lb{rem51}
  By {\rm\x{f9}}, using {\rm\x{f4}} and {\rm\x{f5}}, we obtain
  \beaa S_{1}P(t_{0})S_{1}^{-1}\EQ(Id -P(t_{0})+\hat{P}_{+}(t_{0}))P(t_{0})(Id +P(t_{0})-\hat{P}_{+}(t_{0}))\\
  \EQ \hat{P}_{+}(t_{0}) =U(t_{0},t_{0}).
  \eeaa
Thus it follows from {\rm\x{c14}} that
\bea\lb{f27} \hat{P}_{+}(t)=\hat{\Phi}(t,t_{0})U(t_{0},t_{0})\hat{\Phi}(t_{0},t)
=\hat{\Phi}(t,t_{0})S_{1}P(t_{0})S_{1}^{-1}\hat{\Phi}(t_{0},t).
\eea
Meanwhile,  by {\rm\x{f14}}, using {\rm\x{f6}} and {\rm\x{f8}}, we obtain
  \beaa S_{2}Q(t_{0})S_{2}^{-1}\EQ(Id -Q(t_{0})+\hat{Q}_{-}(t_{0}))Q(t_{0})(Id +Q(t_{0})-\hat{Q}_{-}(t_{0}))\\
  \EQ \hat{Q}_{-}(t_{0}) =V(t_{0},t_{0}).
  \eeaa
  Thus it follows from {\rm\x{d11}} that
\beaa \hat{Q}_{-}(t)=\hat{\Phi}(t,t_{0})V(t_{0},t_{0})\hat{\Phi}(t_{0},t)
=\hat{\Phi}(t,t_{0})S_{2}Q(t_{0})S_{2}^{-1}\hat{\Phi}(t_{0},t),
\eeaa
and consequently,
\bea\lb{f28} \hat{P}_{-}(t)=\hat{\Phi}(t,t_{0})V(t_{0},t_{0})\hat{\Phi}(t_{0},t)
=\hat{\Phi}(t,t_{0})S_{2}P(t_{0})S_{2}^{-1}\hat{\Phi}(t_{0},t).
\eea
By {\rm\x{f21}}, {\rm\x{f27}} and {\rm\x{f28}}, we know that linear operators $\hat{P}_{+}(t)$, $\hat{P}_{-}(t)$ and $\tilde{P}(t)$, defined on $[t_{0},+\oo)$, $(-\oo,t_{0}]$ and $\R$ respectively, are actually obtained under the same rules.
\end{remark}

\begin{remark}
  Throughout this paper we choose any fixed $t_{0} \in \R$ instead of $0\in \R$, which is a little different from the one given in uniform exponential dichotomy (see e.g., {\rm \cite{pop-06}}), where the initial point $0$ is used for simplicity, and there is no substantial difference  in inequalities thus obtained.
  However, here we have to choose general term $t_{0}$ instead of $0$ since the nonuniform item will vanish at time $0$, and hence there is a significant difference in some calculations.
 \end{remark}

\section{Example}
\setcounter{equation}{0} \noindent
In   what   follows   we  use  an   example   to  demonstrate our results.
The following   example
 shows that there exists a linear SDE which admits an
 NMS-ED  but not  uniform.

\begin{example}\lb{exp51} Let $a>b>0$ be real parameters. Then the following linear SDE
 \begin{equation}\lb{xx3}
\left\{ \begin{array}{ll}
du & =(-a-bt\sin t)u(t)dt+\sqrt{2b\cos t}\exp(-at+bt\cos t)d\omega(t)\\
dv & =(a+bt\sin t)v(t)dt-\sqrt{2b\cos t}\exp(at-bt\cos t)d\omega(t)
\end{array} \right.
\end{equation}
with the initial condition $u(0)=v(0)=1$
admits an NMS-ED  that is not a uniform MS-ED.
 \end{example}

\noindent\prf{Let
\[\Phi(t)=\left(
  \begin{array}{cc}
    U(t) & 0 \\
    0 & V(t) \\
  \end{array}
\right)\]
be a fundamental matrix  solution of \x{xx3}. Thus we have $u(t)=U(t)U^{-1}(s)u(s)$ and $v(t)=V(t)V^{-1}(s)v(s)$. In addition, it is easy to verify that
\[\left(
    \begin{array}{cc}
      \exp\left(-at+bt\cos t-b\sin t\right) & 0 \\
      0 & \exp\left(at-bt\cos t+b\sin t\right) \\
    \end{array}
  \right)
\]
is a fundamental matrix solution of
\[
\left\{ \begin{array}{ll}
du & =(-a-bt\sin t)u(t)dt,\\
dv & =(a+bt\sin t)v(t)dt.
\end{array} \right.
\]
Hence, by \cite[p. 97]{ev}, the solution of \x{xx3} is given by
\[
\left\{ \begin{array}{ll}
u(t) & =\exp\left(-at+bt\cos t-b\sin t\right)\left(1+\sqrt{2b}\int_{0}^{t}e^{b\sin s}\sqrt{\cos s}d\omega(s)\right),\\
v(t) & =\exp\left(at-bt\cos t+b\sin t\right)\left(1-\sqrt{2b}\int_{0}^{t}e^{-b\sin s}\sqrt{\cos s}d\omega(s)\right),
\end{array} \right.
\]
since $u(0)=v(0)=1$. Therefore, \beaa \E\|u(t)\|^{2}\EQ \exp\left(-2at+2bt\cos t-2b\sin t\right)\left(1+2b\int_{0}^{t}e^{2b\sin s}\cos sds\right)\\
\EQ \exp\left(-2at+2bt\cos t\right).\eeaa
Thus, one can obtain
 \[\E\|U(t)U^{-1}(s)\|^{2}=\frac{\E\|u(t)\|^{2}}{\E\|u(s)\|^{2}}=e^{-2a(t-s)+2b(t\cos t-s\cos s)}\] since $\E\|u(s)\|^{2}>0$. It is easy to see that
\[\E\|U(t)U^{-1}(s)\|^{2}=e^{(-2a+2b)(t-s)+2bt(\cos t-1)-2bs(\cos s-1)},\] and thus \be\lb{xx4}\E\|U(t)U^{-1}(s)\|^{2}\le e^{(-2a+2b)(t-s)+2b s}, \quad \forall~ (t,s)\in I^{2}_{\geq}.\ee
Furthermore, if $t =4k\pi$ and $s =3k\pi$
with $k\in \N$, then
\be\lb{xx5}\E\|U(t)U^{-1}(s)\|^{2}= e^{(-2a+2b)(t-s)+2b s}, \quad \forall~ (t,s)\in I^{2}_{\geq}.\ee
Similarly, one can prove that
\be\lb{xx6}\E\|V(t)V^{-1}(s)\|^{2}\le e^{(-2a+2b)(s-t)+2b s}, \quad \forall~ (t,s)\in I^{2}_{\leq},\ee
and
\be\lb{xx7}\E\|V(t)V^{-1}(s)\|^{2}= e^{(-2a+2b)(s-t)+2b s}, \quad \forall~ (t,s)\in I^{2}_{\leq} \ee
if $t =4k\pi$ and $s =3k\pi$
with $k\in \N$.
Thus,  \x{xx3} admits an
 NMS-ED. By \x{xx5} and/or \x{xx7}, the exponential $e^{2b s}$ in \x{xx4}
 and/or \x{xx6} cannot be removed. This shows
that the NMS-ED is not uniform.
\hspace{\stretch{1}}$\Box$
}

\begin{remark}
  The SDE {\rm\x{xx3}}  in Example {\rm\ref{exp51}} admitting an
   NMS-ED is linear in the narrow sense. Following the same idea and method in \cite{zc}, one can establish a general linear SDE,
    which admits an  NMS-ED. For example, let $a>b>0$ be real parameters,
    one can prove the following linear SDE
 \beaa
\left\{ \begin{array}{ll}
du & =(-a-bt\sin t)u(t)dt+u(t)d\omega(t)\\
dv & =(a+bt\sin t)v(t)dt+v(t)d\omega(t)
\end{array} \right.
\eeaa
with the initial condition $u(0)=v(0)=1$
admiting an
   NMS-ED that is not a uniform MS-ED.
\end{remark}


\begin{thebibliography}{99}

%

\bibitem{bp02} L. Barreira, Ya. Pesin, Lyapunov exponents and smooth ergodic theory, University Lecture Series 23, Amer. Math. Soc., 2002.


\bibitem{bp07} L. Barreira, Ya. Pesin,   Nonuniform  Hyperbolicity,  Encycl.  Math.  Appl.,   vol. 115,   Cambridge  Uni-versity  Press,   2007.

\bibitem{bcv11} L.   Barreira,  J.  Chu,  C.  Valls, Robustness of nonuniform dichotomies with different growth rates, S\~{a}o Paulo J. Math. Sci, 5  (2011), 203-231.

\bibitem{bcv13} L.   Barreira,  J.  Chu,  C.  Valls, Lyapunov Functions for General Nonuniform Dichotomies, Milan J. Math, 81 (2013), 153-169.

\bibitem{bsv09} L.   Barreira,  C.  Silva,  C.  Valls,  Nonuniform  behavior   and   robustness,  J.   Differential  Equations,  246
(2009),   3579-3608.


\bibitem{bv06} L.  Barreira, C. Valls, Stable  manifolds for nonautonomous equations without exponential dichotomy,
J.   Differential  Equations,  221   (2006),   58-90.


\bibitem{bv08}L. Barreira, C. Valls, Robustness of nonuniform exponential dichotomies in Banach spaces, J.   Differential  Equations, 244 (2008), 2407-2447.


\bibitem{chi} C. Chicone, Yu. Latushkin, Evolution Semigroups in Dynamical
Systems and Differential Equations, Mathematical Surveys and Monographs 70, Amer. Math. Soc. 1999.

\bibitem{cl-94} S. N. Chow, H. Leiva,  Dynamical   spectrum  for  time  dependent   linear   systems   in   Banach  spaces,
Jpn.  J.   Ind.   Appl.  Math.,  11  (1994),   379-415.



\bibitem{cl} S. N. Chow, H. Leiva, Existence and roughness of the exponential
dichotomy for skew-product semiflows in Banach spaces, J.   Differential  Equations, 120 (1995), 429-477.



\bibitem{cof-71}C. V.  Coffman,  J.J.   Sch\"{a}ffer,   Linear   differential   equations  with  delays:  admissibility   and   conditional  exponential  stability,
J.   Differential  Equations,  9  (1971),   521-535.

\bibitem{cop} W. A. Coppel, Dichotomy in stability theory, Lecture Notes in Mathematics, Vol. 629,
Springer-Verlag, New York/Berlin, 1978.

\bibitem{dk-74} J. L. Dalec'ki\u\i, M. G. Kre\u\i n,  Stability   of   Differential  Equations  in   Banach  Space, Amer. Math. Soc., Providence, R.I., 1974.

\bibitem{drk}T.S. Doan, M. Rasmussen, P.E. Kloeden, The mean-square dichotomy spectrum and a bifurcation to a mean-square attractor, Discrete Contin. Dyn. Syst. Ser. B, 20 (2015), 875-887.

 \bibitem{en}K. J. Engel, R. Nagel, One-parameter semigroups for linear evolution equations, Springer-Verlag, 2000.

\bibitem{ev} L. C. Evans, An introduction to stochastic differential
equations, Amer. Math. Soc. 2012.



\bibitem{fl} M. Fu, Z. Liu, Square-mean almost automorphic solutions for some
stochastic differential equations, Proc. Amer. Math. Soc., 138
(2010), 3689-3701.


\bibitem{hale} J. K. Hale, Ordinary  differential  equations,
Wiley-Interscience, New York, 1969.

\bibitem{hl-86} J. Hale, X. B. Lin, Heteroclinic orbits for retarded functional differential equations, J. Differential Equations, 65 (1986), 175-202.

\bibitem{hen-81}D. Henry, Geometric Theory of Semilinear Parabolic Equations, Lecture Notes in Math. 840, Springer-Verlag, Berlin, 1981.

\bibitem{hig} D. J. Higham, Mean-square and asymptotic stability of
the stochastic theta method, SIAM J. Numerical Anal., 38
(2000), 753-769.

\bibitem{hms} D. J. Higham, X. Mao, A. M. Stuart, Exponential
mean-square stability of numerical solutions to stochastic
differential equations, LMS J. Comput. Math., 6 (2003), 297-313.

\bibitem{hmy} D. J. Higham, X. Mao, C. G. Yuan, Preserving exponential
mean-square stability in the simulation of hybrid stochastic differential
equations, Numer. Math., 108 (2007), 295-325.

\bibitem{huy-06}N.T. Huy, Exponential dichotomy  of  evolution equations and  admissibility  of  function  spaces  on a half-line, J.  Funct.  Anal.,
235   (2006),   330-354.

\bibitem{il-01}P. Imkeller, C. Lederer, On the cohomology of flows of stochastic and random differential equations, Prob. Theory
Related Fields, 120 (2001), 209-235.

\bibitem{jw-01}N.  Ju,  S.   Wiggins,  On  roughness   of   exponential  dichotomy,  J.   Math.  Anal.  Appl.,  262   (2001),   39-49.

\bibitem{kl} P. E. Kloeden, T. Lorenz, Mean-square random dynamical
systems, J. Differential Equations, 253 (2012), 1422-1438.

\bibitem{lad} A. G. Ladde, G. S. Ladde, An introduction to differential equations, Volume 2.
Stochastic Modeling, Methods and Analysis, World Scientific Publishing Co. 2013.

\bibitem{lmr}Y.  Latushkin,   S.   Montgomery-Smith,  T.  Randolph,   Evolutionary semigroups and dichotomy of linear skew-product flows on locally compact spaces with Banach fibers, J.   Differential  Equations, 125   (1996),
73-116.

\bibitem{lrs}Y.  Latushkin,   T.  Randolph,   R.  Schnaubelt,  Exponential  dichotomy   and   mild  solutions  of   nonautonomous  equations  in
Banach  spaces,  J.   Dynam.  Differential  Equations,  10  (1998),   489-510.

\bibitem{lin} X. B.  Lin,   Exponential  dichotomies  and   homoclinic   orbits   in   functional   differential   equations,   J.   Differential  Equations,  63
(1986),   227-254.

\bibitem{lin-94} X. B.  Lin,   Exponential dichotomies in intermediate spaces with applications to a diffusively perturbed predator-prey model, J. Differential Equations, 108 (1994), 36-63.

\bibitem{ls} Z. Liu, K. Sun, Almost automorphic solutions for stochastic
differential equations driven by L\'{e}vy noise, J. Funct. Anal., 226
(2014), 1115-1149.

\bibitem{liz-92} M. Lizana, Exponential dichotomy for singularly perturbed linear functional differential equation with small delays, Appl. Anal.,  47 (1992), 213-225.

\bibitem{mao} X. Mao, stochastic differential equations and
applications, Horwood, Chichester, 1997.

\bibitem{ms} J. Massera, J. Sch\"{a}ffer, Linear differential equations and
functional analysis  I,  Ann. of Math., 67 (1958), 517-573.


\bibitem{np} R.  Naulin,  M.  Pinto, Roughness of $(h,k)$-dichotomies,  J.   Differential  Equations,  118   (1995),   20-35.

\bibitem{np-97} R.  Naulin,  M.  Pinto,   Stability   of   discrete  dichotomies  for  linear   difference  systems,   J.   Difference  Equ.  Appl.,  3  (1997),
101-123.

\bibitem{pa-84}K. J. Palmer, Exponential dichotomies and transversal homoclinic points, J. Differential Equations, 55 (1984), 225-256.

\bibitem{pa-88}K. J. Palmer, Exponential dichotomies and Fredholm operators, Proc. Amer. Math. Soc., 104 (1988), 149-156.

\bibitem{per} O. Perron, Die Stabilit\"{a}tsfrage bei Differentialgleichungen, Math. Z., 32 (1930), 703-728.

\bibitem{per-71} G.  Pecelli,  Dichotomies  for  linear   functional-differential   equations,   J.   Differential  Equations,  9  (1971),   555-579.

\bibitem{ps-99}V. Pliss, G. Sell,  Robustness of  exponential dichotomies in  infinite-dimensional dynamical  systems,  J.  Dynam. Differential
Equations,  11  (1999),   471-513.

\bibitem{pop-06}L. H. Popescu, Exponential dichotomy roughness on Banach spaces, J. Math. Anal. Appl., 314 (2006), 436-454.

\bibitem{pop-09}L. H. Popescu, Exponential dichotomy roughness and structural stability for evolution families without bounded growth and decay. Nonlinear Anal., 71 (2009), 935-947.

\bibitem{pp-04}P. Preda,  A.  Pogan,  C.  Preda, On {$(a,b)$}-dichotomy for evolutionary processes on a half-line, Glasg.
Math.  J.,   46  (2004),   217-225.


\bibitem{pp-06}P. Preda,  A.  Pogan,  C.  Preda,  Sch\"{a}ffer   spaces   and   exponential  dichotomy   for  evolutionary   processes,  J.   Differential
Equations,  230   (2006),   378-391.



\bibitem{rr-95}H. M.  Rodrigues,   J. G.   Ruas-Filho,  Evolution  equations:   dichotomies  and   the   Fredholm  alternative   for  bounded   solutions,
J.   Differential  Equations,  119   (1995),   263-283.

\bibitem{ss} R. Sacker, G. Sell, Dichotomies for linear evolutionary equations
in Banach spaces, J.   Differential  Equations, 113 (1994), 17-67.

 \bibitem{ss1}  R. Sacker, G. Sell, Existence of dichotomies and invariant
splitting for linear differential systems I [II, III], J.   Differential  Equations, 15 (1974), 429-458 [22 (1976), 478-496,
497-522].

\bibitem{sk-01} O. M. Stanzhyts'kyi, Exponential dichotomy and mean square bounded solutions of linear stochastic Ito systems, Nonlinear Oscil., 4 (2001), 389-398.

\bibitem{sk-06} O. M. Stanzhyts'kyi, A. P. Krenevych, Investigation of the exponential dichotomy of linear stochastic It\^{o} systems with random initial data by means of quadratic forms, Ukrainian Math. J., 58 (2006), 619-629.

\bibitem{st} D. Stoica, Uniform exponential dichotomy of stochastic
cocycles, Stochastic Process. Appl., 120 (2010), 1920-1928.

\bibitem{zlz} L. Zhou, K. Lu, W. Zhang, Roughness of tempered exponential dichotomies for infinite-dimensional random difference equations, J.   Differential  Equations, 254 (2013), 4024-4046.

\bibitem{zlz-17} L.  Zhou, K. Lu,  W. Zhang, Equivalences between nonuniform exponential dichotomy and admissibility, J.   Differential  Equations, 262 (2017), 682-747.

\bibitem{zz-16} L.  Zhou,   W. Zhang,  Admissibility  and   roughness   of   nonuniform  exponential  dichotomies  for  difference  equations,
J. Funct.  Anal.,  271   (2016),   1087-1129.



\bibitem{zc1}H. Zhu, J. Chu, Mean-square exponential dichotomy of numerical solutions to
stochastic differential equations, J. Appl. Anal. Comput., 6 (2016), 463-478.

\bibitem{zc2}H. Zhu, Y. Jiang, Robustness of mean-square exponential dichotomies for linear stochastic equations, Electron. J. Differential Equations, 123 (2017), 1-13.

\bibitem{zcz}H. Zhu, J. Chu, W. Zhang, Mean-square Almost automorphic solutions for
stochastic differential equations with hyperbolicity, Discrete Contin. Dyn. Syst., 38(4) (2018), 1935-1953.

\bibitem{zc}H. Zhu,   L. Chen, Nonuniform exponential dichotomies in mean square and second-moment Lyapunov exponent, preprint.


\end{thebibliography}
\end{document}